\documentclass[]{article}
\usepackage{latexsym}
\usepackage{mathptmx}
\usepackage{graphicx}
\usepackage{blindtext}
\usepackage{tabularx} 
\usepackage{rotating}
\usepackage{setspace}
\usepackage{amsmath}
\usepackage{float}
\usepackage{mathptmx}
\usepackage{breakcites}
\usepackage{amsbsy}
\usepackage{enumerate}
\usepackage{amssymb}
\usepackage{color}
\usepackage{graphicx}
\usepackage{tabularx} 
\usepackage{placeins}
\usepackage{amsthm}
\usepackage{mathptmx}
\usepackage{csquotes}
\usepackage[titletoc,title]{appendix}
\MakeOuterQuote{"}

\newtheorem{theorem}{Theorem}
\newtheorem{assumption}{Assumption}
\newtheorem{corollary}[theorem]{Corollary}
\newtheorem{lemma}{Lemma}

\DeclareMathOperator*{\argmin}{arg\,min}
\DeclareMathOperator*{\argmax}{arg\,max}

\title{Hesitant  Adaptive Search with Estimation and Quantile Adaptive Search for Global Optimization with Noise}
\author{David D. Linz, Zelda B. Zabinsky\thanks{zelda@uw.edu} \\ Department of Industrial and Systems Engineering \\
University of Washington, Seattle, Washington 98195}

\begin{document}

\maketitle


\begin{abstract}
Adaptive random search approaches have been shown to be effective for global optimization problems, where under certain conditions, the expected performance time increases only linearly with dimension. However, previous analyses assume that the objective function can be observed directly. We consider the case where the  objective function must be estimated, often using a noisy function, as in simulation. We present a finite-time analysis of algorithm performance that combines estimation with  a sampling distribution. 
We present a framework called Hesitant Adaptive Search with Estimation, and derive an upper bound on function evaluations that is cubic in dimension, under certain conditions.  We extend the framework to  Quantile Adaptive Search with Estimation, which focuses sampling points from a series of nested quantile level sets. The analyses suggest that computational effort is better expended on sampling improving points than refining estimates of objective function values during the progress of an adaptive search algorithm.
\end{abstract}

\section{Introduction}
\label{sec:intro}

Adaptive random search algorithms for global optimization are often characterized by how they generate points within a feasible region, i.e., their sampling distribution. Many algorithms attempt to iteratively improve the  sampling distribution to focus on promising regions, based on observations   \cite{boender1995stochastic,FuHandbook,HuJiaqiao2012ASoS, locatelli2013global,LocatelliSchoen2021,pardalos2000recent,raphael2003direct,RaphaelSmith,CrossEntropy,Zabinsky2003}.  
An additional complication is the need to estimate the objective function value  as the algorithm progresses.

This paper provides a finite-time analysis of a class of adaptive random search algorithms applied  to problems that require {\it estimation} to account for noise, that is, problems where the objective function cannot be evaluated directly but must be estimated. We describe an extension to Hesitant Adaptive Search, we call Hesitant Adaptive Search with Estimation (HAS-E), and embed a confidence interval on the estimate of the current objective function value into the algorithm.  The analysis relates the number of replications used in the estimation of the objective function to the overall performance of the algorithm.

In contrast to asymptotic convergence, a finite-time analysis of performance is available for a class of adaptive random search algorithms when the objective function can be evaluated directly through an oracle (i.e., a black-box function).  Specifically, prior research has derived a finite-time analysis for 
 Pure Adaptive Search ({PAS}), Hesitant Adaptive Search ({HAS}), Backtracking Adaptive Search ({BAS}), and Annealing Adaptive Search ({AAS}) \cite{BAS2007,Bulger1998,romeijn1994simulated,Shen2005,BAS2006,Wood2001,Zabinsky2003,Zabinsky2010,Zabinsky1992,Zabinsky1995}. 
 These algorithms provide a framework for analysis, and are not intended to be implemented directly. However, the analyses shed light on the role of the sampling distribution and probability of generating improving points on performance.
 Under certain conditions, the expected number of function evaluations required to sample below a specified objective function value increases only linearly in dimension when optimizing a function without noise.
We address the question of how estimation of a noisy function impacts the performance.

We first analyze the performance of HAS-E, where the sampling distribution focuses on nested level sets while allowing for `hesitation.'  We also introduce a new adaptive random search algorithm, called Quantile Adaptive Search with Estimation (QAS-E), which samples over the entire domain but parametrically modifies the sampling distribution based on a sequence of quantiles.  The motivation for the analysis of QAS-E is to provide finite-time analyses that can be adapted for use in  adaptive random search algorithms that use quantiles in the adaptive mechanics \cite{HoOrdinal2000,Ho:OrdinalBook,JiangHuPeng2022quantilebased,zabinsky2019PBnB}. 

The main result of this paper is in
Theorem \ref{thm:bound_hase}, which provides an upper bound on the expected number of function evaluations (including replications) required to first obtain a value within a target $\epsilon$ of the global minimum. 
This is used to show, in
 Corollary \ref{thm:hase_cor}, that under certain conditions the expected number of
function evaluations (including replications) to obtain a value less than $\epsilon$ above the minimum is bounded by a cubic function of the domain dimension. We then use the analysis of HAS-E to derive analogous bounds for QAS-E, which resembles quantile-based algorithms in practice.

\section{Preliminaries}
\label{sec:background}

Consider an optimization problem,

\noindent
\begin{equation} \label{eq:main_prob_zero} 
\min_{x \in S} f(x)  \tag{\textit{P}}
\end{equation}

\noindent
where $S$ is a closed and bounded subset of $\mathbb{R}^n$, $x \in S \subset \mathbb{R}^n$, and $f: \mathbb{R}^n \rightarrow \mathbb{R}$. Denote the minimum value and an optimal point in the domain, 
respectively, as:
\begin{equation}\label{eq:max_def}  y_* = \min_{x \in S} f(x) \text{ and } x_* = \argmin_{x \in S}  f(x) . \end{equation}
\noindent
Similarly, denote the maximum value and point as:
\noindent
\begin{equation}\label{eq:min_def}  y^* = \max_{x \in S } f(x) \text{ and } x^* = \argmax_{x \in S }  f(x) . \end{equation}
\noindent
Furthermore, we define the diameter $d$ of $S$ as the greatest distance between any two points in $S$. 	

We are particularly interested in the situation where the objective function cannot be evaluated directly but is expressed as  
\begin{equation}\label{eq:noisyf}
f(x)  = E [ g(x, \chi) ] 
\end{equation}
where $g(x, \chi)$ is a \enquote{noisy} function of $x\in S$ and a random variable $\chi$. Often, $g(x, \chi)$ is evaluated using a discrete-event simulation and is replicated a number of times at each point $x$, taking the sample mean as an estimate of $f(x)$.  In this paper, we seek to relate the number of replications on an iteration to the overall performance of an adaptive random search algorithm.

The finite-time analysis of Pure Adaptive Search establishes that, under certain conditions, the expected number of iterations (i.e., function evaluations) until PAS achieves a value close to the minimum increases only linearly in terms of the dimension of the domain 
\cite{Zabinsky2003, Zabinsky1992, Zabinsky1995}.  
However, the requirement that PAS improves at every iteration makes it difficult to implement practically. Hesitant Adaptive Search generalizes PAS by relaxing this requirement and allowing hesitation, thereby extending the class of algorithms it represents \cite{Bulger1998,Wood2001}.
%
We summarize HAS and a key result. 

HAS is defined by a sampling distribution $\zeta$ with support on $S$ and a bettering probability $b(y)$, $0<b(y)\leq1$, defined for $y_* < y \leq y^*$.  On any iteration with objective function value $y$, HAS generates an improving point with probability $b(y)$ by drawing from the normalized restriction of $\zeta$ on the improving level set. The probability of hesitation is $1 - b(y)$, where the current point does not change.   HAS is defined as follows. 

\vspace{3mm} 

\begin{samepage}
\noindent
\textbf{Hesitant Adaptive Search (HAS),  cf. \cite{Bulger1998}   }
\begin{itemize} 
	\item \textbf{Step 0:} Sample $X_0$ in $S$ according to the probability distribution $\zeta$ on $S$. Set $\bar{Y}_0 =f(X_0)$. Set $k = 0$. 
	\item \textbf{Step 1:} Generate $X_{k+1}$ from the normalized restriction of $\zeta$ on the improving set $S_{k} = \{ x \in S: f(x) < \bar{Y}_k \} $ with probability $b(\bar{Y}_k)$, and set $\bar{Y}_{k+1} = f(X_{k+1})$. Otherwise, set $X_{k+1} = X_k$ and $\bar{Y}_{k+1} = \bar{Y}_{k}$.
	\item \textbf{Step 2:} If a stopping criterion is met, stop. Otherwise, increment $k$ and return to Step 1.
\end{itemize}
\end{samepage}
Note that $\bar{Y}_{k}$ is 
non-increasing, and is decreasing on any iteration with probability $b(\bar{Y}_{k})$.

The finite-time analysis of HAS in \cite{Bulger1998} provides a closed form expression for a bound on the expected number of iterations until reaching a specified $\epsilon > 0$ above the minimum function value, denoted $E[N(y_* +\epsilon)]$, as,
\begin{equation} \label{eq:has_bound_origin} E[N(y_* +\epsilon)] \leq 1 + \int_{y_* + \epsilon}^{\infty} \frac{d\rho (t)}{b(t) \cdot p(t)}  \end{equation}
\noindent 
where $\rho(y) = \zeta(f^{-1}([ -\infty, y ]))$, and $p(y) = \rho((-\infty, y ]) $.  A complete characterization of HAS for problems with mixed continuous-integer variables is in \cite{Wood2001}. 

%

The HAS analysis provides insight into the relationship between the bettering probability and performance. 
However, HAS is still difficult to implement because it is impractical to draw from the normalized restriction of $\zeta$ on the improving level set. 
Another way to define an adaptive random search algorithm is to always sample from the entire set $S$, and iteratively update a parameter controlling the sampling distribution.

Annealing Adaptive Search is an abstraction of simulation annealing, and it always samples from a Boltzmann distribution on the entire set $S$. The temperature parameter for the Bolzmann distribution is iteratively decreased to control the  update of the sampling distribution. The analysis in \cite{Shen2005,Shen2007} establishes stochastic dominance between
AAS and a special case of HAS, making use of the finite-time analysis of HAS.
 The analysis of AAS was  used to derive an analytical cooling schedule for simulated annealing algorithms \cite{Shen2007}. 

The analyses of PAS, HAS, and AAS provide insight into the performance of adaptive random search algorithms, however, they assume the objective function $f(x)$ can be evaluated exactly.  As random search algorithms are being applied broadly to functions that require estimation, a major question is how estimation impacts performance.

%
\section{Hesitant Adaptive Random Search with Estimation (HAS-E)}
\label{sec:hase}

\vspace{2mm}
\noindent


Given that the value of $f(x)$, for $x \in S$, cannot be directly observed, we consider estimating the value by performing a certain number of independent replications and taking the sample mean. Suppose $g(x, \chi_r)$ is evaluated at a point $x$ for $R$ replications, $r =1, \ldots, R$. The sample mean estimate (dropping the $x$ for notational convenience) is,
\begin{equation}\label{eq:basic_estimation} \hat{y}^{est} = \frac{\sum_{r=1}^{R} g(x, \chi_r)}{R}. \end{equation}
\noindent
We assume that $\hat{y}^{est} \sim N(f(x), \frac{\sigma}{\sqrt{R}})$, where $\sigma^2 = Var(g(x, \chi))$, and $\sigma$ is known. 
A standard probability bound is given by
\begin{equation} \label{eq:yestprob}
P\left(f(x) - \frac{\sigma \cdot z_{\alpha/2}}{\sqrt{R}} \leq \hat{y}^{est} \leq f(x)  + \frac{\sigma \cdot z_{\alpha/2}}{\sqrt{R}}  \right) \geq 1 - \alpha
\end{equation}
for $ 0 \leq \alpha \leq 1$ and where $z_{\alpha/2}$ is the standard normal value at $\alpha/2$. 

%
%

We are interested in an upper bound of the estimate and let $\hat{y}^{high}$ be the upper confidence interval value, given by
%
\begin{equation} \label{eq:basic_estimate}
\begin{split}
& \hat{y}^{high} = \hat{y}^{est} + \frac{\sigma \cdot z_{\alpha/2}}{\sqrt{R}}.  \\ 
\end{split}
\end{equation}
Since we want to know how far $\hat{y}^{high}$ is from the true value $f(x)$, we  note that $ \hat{y}^{high}  \sim N\left(f(x) + \frac{\sigma \cdot z_{\alpha/2}}{\sqrt{R}} , \frac{\sigma}{\sqrt{R}}\right)$ and,  also from (\ref{eq:yestprob}), we have,

\begin{equation}\label{eq:conf} P\left( f(x) \leq \hat{y}^{high} \leq f(x) + 2\frac{\sigma \cdot z_{\alpha/2}}{\sqrt{R}}  \right) \geq 1 - \alpha .\end{equation}
\noindent 

HAS-E uses the estimate $\hat{y}^{high}$ to focus sampling on regions that are likely to be improving. In contrast to HAS that samples in the improving level set with a bettering probability, HAS-E samples in the level set associated with $\hat{y}^{high}$ with a bettering probability.

On the $k$th iteration of HAS-E, the sampled point $x_k \in S$ is evaluated with $R_k$ independent replications of $g(x_k, \chi_r)$ for $r=1, \ldots, R_k$, and then the estimate  $\hat{y}_k^{est}$ is calculated as in (\ref{eq:basic_estimation}) and the upper bound $\hat{y}_k^{high}$ as in (\ref{eq:basic_estimate}).   We let
 $y_k = f(x_k)$  be the true objective function value at $x_k$, and the improving level set $S_{y_k}$ be\begin{equation} \label{eq:level_set_def} S_{y_k} = \{x \in S: f(x) < y_k  \}. \end{equation}
Similarly, we let
 \begin{equation} \label{eq:level_set_upperbound} S_{\hat{y}_k^{high}} = \{x \in S: f(x) < \hat{y}_k^{high}  \} \end{equation}
be the level set associated with the upper confidence interval bound and note that 
$$P( S_{y_k}\subset S_{\hat{y}_k^{high}}) \geq 1-\alpha$$
from \eqref{eq:conf}.

\begin{figure}[!ht]
	\centering
	\includegraphics[height=3.1in]{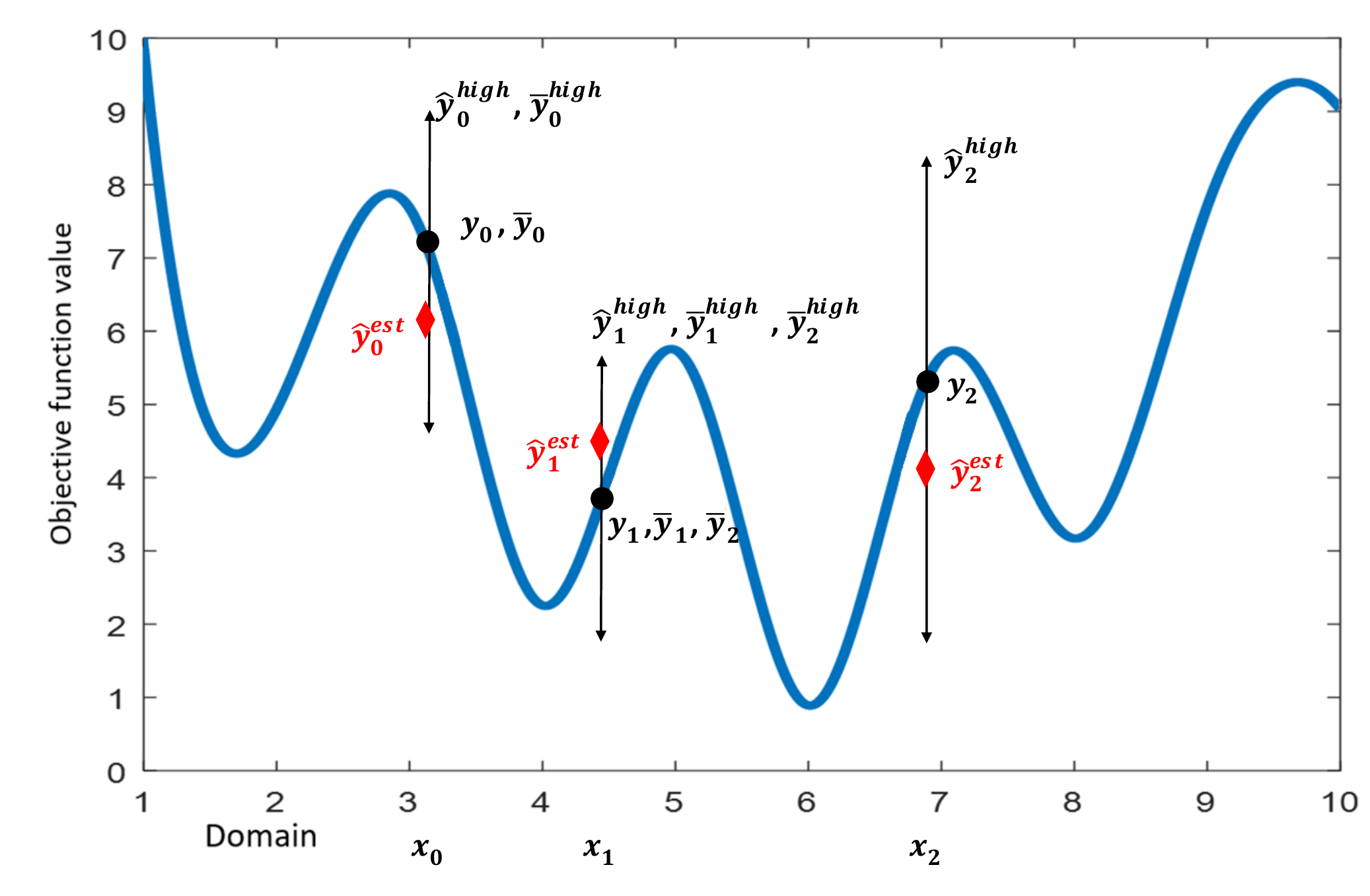} 
	\caption{An illustration of three HAS-E  iterations ($k= 0,1,2$) on a one-dimensional problem. The sampled points are labeled $x_0, x_{1},$ and $x_{2}$, with true values $y_0, y_{1},$ and $y_{2}$. The estimated values, $\hat{y}^{est}_k$, and the upper confidence values, $\hat{y}_k^{high}$, are shown for $k= 0,1,2$. The best true values, $\bar{y}_0, \bar{y}_{1},$ and $\bar{y}_{2}$, and their associated upper confidence values $\bar{y}_0^{high}, \bar{y}_{1}^{high},$ and $\bar{y}_2^{high}$ are also illustrated. \\ \vspace{2mm}
	}\label{fig:HASE_example1}
\end{figure}

We are interested in sampling a point in a target level set $S_{y_*+\epsilon}$ for some $\epsilon>0$.  The general approach of HAS-E is to sample from the normalized restriction of $\zeta$ on $S_{\bar{y}_k^{high}}$
with bettering probability $b(\bar{y}_k^{high})$, and hesitate (remain at the same point) with probability $1-b(\bar{y}_k^{high})$, where $\bar{y}_k^{high}$ is the estimated upper confidence bound on the $k$th iteration.   Later, for analysis purposes, we let $\gamma$ be a minimum bettering probability such that $b(y)\geq \gamma$ for $y_*+\epsilon \leq y \leq y^*$. See Figure~\ref{fig:HASE_example1} for an illustration of three iterations of HAS-E. 

HAS-E requires input parameters $\alpha$ and  $\sigma$, along with a sampling distribution $\zeta$ with support on the entire domain $S$, and the bettering probability $b(y)$.  It also requires a sequence of the number of replications on iteration $k$, i.e., $\{R_k, k=0, 1,\ldots \}$, which is discussed later.

As with HAS, the HAS-E algorithm is a framework for analysis, and not intended to be implemented directly. However, the framework allows us to analyze the algorithm's performance on any iteration $k$.

\vspace{4mm}    
\noindent
\textbf{Hesitant Adaptive Search with Estimation (HAS-E)}
%
\begin{itemize} 
	\item \textbf{Step 0:} Sample $X_0$ in $S$ according to the probability distribution 
	$\zeta$ 
	on $S$. Conduct $R_0$ independent replications of the function at the initial selected point, i.e., $g(X_0, \chi_r) $ for $r = 1, \ldots, R_0$. Estimate the value $\hat{y}_0^{high} $ as in (\ref{eq:basic_estimate}) and set $\bar{y}_0^{high} = \hat{y}_0^{high}$.  
	Set $\bar{Y}_0=f(X_0)$.	Set $k = 0$. 
	
	\item \textbf{Step 1:} Generate $X_{k+1}$ from  the normalized restriction of 
	$\zeta$ 
	on the set $S_{\bar{y}_k^{high}}$ with bettering probability $b(\bar{y}_k^{high})$, and estimate $\hat{y}_{k+1}^{high}$ as in (\ref{eq:basic_estimate}) with $R_k$ independent replications of $g(X_{k+1}, \chi_r) $ for $r = 1, \ldots, R_k$.
	Otherwise (with probability $1-b(\bar{y}_k^{high})$), set $X_{k+1} = X_k$ 
	and $\hat{y}_{k+1}^{high}=\hat{y}_{k}^{high}$. 	
	Then update  
$$	\bar{Y}_{k+1} =
\begin{cases}
f(X_{k+1}) & \text{if}\ f(X_{k+1}) < \bar{Y}_{k} \\
\bar{Y}_{k} & \text{otherwise}
\end{cases} $$
	and its associated upper confidence bound estimate,
\begin{equation}
	\bar{y}_{k+1}^{high}= \left\{ \begin{array}{ll}
	     \hat{y}_{k+1}^{high} & {\rm if}\  f(X_{k+1}) < \bar{Y}_{k} \\
	     \bar{y}_{k}^{high} & {\rm otherwise.}
	     \end{array} \right.
	     \nonumber
\end{equation}

	%
%
	\item \textbf{Step 2:} If a stopping criterion is met, stop. Otherwise, increment $k$ and return to Step $1$.
\end{itemize}
Note that $\bar{Y}_k$ is  non-increasing, however it is possible for $ \bar{y}_{k}^{high}$ to increase and decrease.

The analysis begins with Theorem \ref{thm:lower_bound_convex}, characterizing the number of replications chosen on each iteration. We next prove in Theorem \ref{thm:hase_stoch_dom} that under certain assumptions about the replications, HAS-E stochastically dominates a special case of HAS without estimation. This allows us to provide (in Theorem \ref{thm:bound_hase})  upper bounds on the expected number of HAS-E iterations and the expected number of function evaluations, including replications, to achieve an optimal solution with function value below a specified threshold. 
Finally, Corollary \ref{thm:hase_cor} provides  bounds on the expected number of HAS-E iterations and  expected number of function evaluations, that, under certain assumptions, increase linearly in dimension and cubic in dimension, respectively.

For purposes of our analysis, we consider a lower bound on the ratio of volumes of the true level set to the level set associated with the upper confidence interval bound $\hat{y}_k^{high}$, i.e., ${\nu(S_{y_k})}/{\nu(S_{\hat{y}_k^{high}})}, $  where $\nu(\cdot)$ is the $n$-dimensional volume of a set. A trivial lower bound for this quantity can be based on the desired accuracy value $\epsilon$,  with $\epsilon> 0$. For $y_k$ and $\hat{y}_k^{high}$ such that  $y_* + \epsilon < y_k \leq \hat{y}_k^{high} \leq y^* $, the following lower bound on the ratio 
$$ \frac{\nu(S_{y_k})}{\nu(S_{\hat{y}_k^{high}})}  \geq  \frac{\nu(S_{y_*+\epsilon}) }{\nu(S)} $$
\noindent 
 holds.

\begin{figure}[!ht]
	\centering
	\includegraphics[height=3.1in]{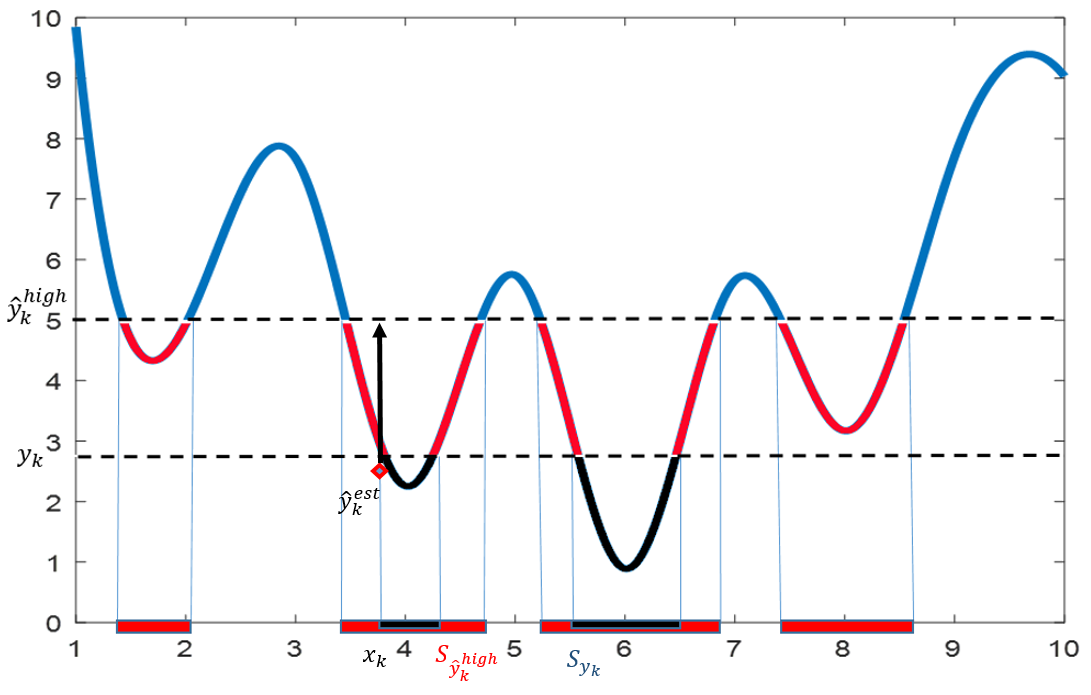}
	\caption{An illustration of the values $y_{k}$ and $\hat{y}_k^{high}$ along with their corresponding level sets 	$S_{y_{k}}$ and $S_{\hat{y}_k^{high}}$ for a one-dimensional problem. The level sets are shown highlighted on the horizontal axis. The ratio between the volumes of the level sets, ${\nu(S_{y_{k}})}/{\nu(S_{\hat{y}_k^{high}})}$, increases and approaches one as the difference between $\hat{y}_k^{high}$ and $y_{k}$ decreases and $\hat{y}_k^{high}$ approaches $y_k$ (which happens with a large number of replications).\\ \vspace{2mm}}
	\label{fig:HASE_example_ratio}
\end{figure}

However, this trivial lower bound may be very small, so we consider another lower bound, and relate it  to  the number of replications used in the estimation.  We let $q$ denote a lower bound on the ratio, with  $0 < q < 1$, 
such that 
$$\frac{\nu(S_{y_k})}{\nu(S_{\hat{y}_k^{high}})}  \geq  q$$
for $y_* + \epsilon < y_k \leq \hat{y}_k^{high} \leq y^* $.
%
Figure~\ref{fig:HASE_example_ratio} illustrates the level sets $S_{{y}_k} $ and $S_{\hat{y}_k^{high} }$. Notice that the ratio of the volumes ${\nu(S_{y_k})}/{\nu(S_{\hat{y}_k^{high}})}$ becomes close to $1$ as the distance between the values $\hat{y}_k^{high} $ and $y_k$ decreases, which typically occurs as the number of replications increases. 
Theorem \ref{thm:lower_bound_convex} provides a bound on the number of replications needed such that the ratio of volumes can be bounded below for a selected value $q$. 

For a function $f(x)$ on domain $S$, we define a quantity $\mathcal{K}_{q}$, for a given $0 < q <1$, which can be viewed as the maximum ratio of the change in objective function to the diameter of $S$
\begin{equation}\label{eq:K_q}
\mathcal{K}_{q} = \frac{\kappa_{q}}{d}  
\end{equation} 
where $\kappa_q$ is the maximum value such that $\nu(S_z) / \nu(S_{z + \kappa_{q}}) >~q$ for any $z$, $y_* < z < y^*$.  The quantity $\mathcal{K}_{q}$ depends on characteristics of the problem.


Furthermore, we define $\mathcal{B}_y$ as the largest ball centered at $x_*$ that can be inscribed inside a level set $S_y$ for  $y_* < y < y^*$, and let  $r_y$ be its radius. Using these two defined concepts, we now relate the number of replications to a selected $q$.

\vspace{ 3mm}
\noindent
\begin{theorem}\label{thm:lower_bound_convex}
	Consider problem (P) and the $k$th iteration of HAS-E with $x_k \in S$, $y_k = f(x_k)$,  $  y_* + \epsilon < y_k \leq y^* $ for $\epsilon > 0 $, and  with $\hat{y}_k^{high}$ estimated with $R $ replications, as in \eqref{eq:basic_estimate}. Also, suppose  that ${y_k} \leq \hat{y}_k^{high} \leq {y_k} + \frac{2\cdot\sigma \cdot z_{\alpha/2}}{\sqrt{R}}$ (which occurs with probability $(1- \alpha)$). For any given value $0 < q < 1$ and associated   $\mathcal{K}_{q}$,  if
	\begin{equation} \label{eq:r_bound}
	R \geq \left( \frac{ \sqrt[n]{q} \cdot 2  \cdot  \sigma \cdot z_{\alpha/2}}{(1 - \sqrt[n]{q}) \cdot r_{y_*+\epsilon} \cdot \mathcal{K}_{q} } \right)^2  
	\end{equation}
	then
	\begin{equation}\label{eq:q_bound} 	\frac{\nu(S_{y_{k}})}{\nu( S_{\hat{y}_k^{high}})} \geq  q.
	\end{equation}

\end{theorem}

\noindent
\textbf{Proof:} See the proof in Appendix~\ref{sec:proofs_hase}. 
\vspace{2mm}

Theorem~\ref{thm:lower_bound_convex} provides an expression for the number of replications needed to achieve a lower bound of $q$. When the value of $q$ is close to one, indicating a relatively tight upper confidence bound, more replications are needed than the number of replications when $q$ has a small value.  The difference becomes more pronounced as the dimension of the problem increases. We ask the question whether a tight estimate is worthwhile for overall performance.

We next develop an upper bound on the expected number of iterations of HAS-E to achieve a function value of $y_*+\epsilon$ or better.  
The analysis of HAS-E proceeds   
in Theorem~\ref{thm:hase_stoch_dom}  by showing that HAS-E stochastically dominates a special case of the HAS algorithm, that we call HAS1.  Then, using HAS1, Theorem~\ref{thm:bound_hase}  provides an upper bound on the expected number of iterations of HAS-E and expected number of function evaluations (including replications) to achieve $y_*+\epsilon$ or better.


The special case HAS1 has a uniform sampling distribution, i.e., $\zeta^{HAS1} \sim {\rm Uniform}$,  and 
the bettering probability is chosen to be constant for all $y$, 
\begin{equation}
b^{HAS1}(y)= \gamma \cdot ( 1 - \alpha) \cdot q  
\label{eq:HAS1betteringprob}
\end{equation} 
where $0<\gamma\leq1$,
 $0<\alpha<1$,  and $0<q<1$.

Let $\bar{Y}_k^{HASE}$ be the best sampled value on the $k$th iteration of the HAS-E algorithm. 
Let $\bar{Y}_k^{HAS1}$ be the best sampled value of HAS1 on the $k$th iteration. 

For the performance analysis of HAS-E in Theorems~\ref{thm:hase_stoch_dom}, \ref{thm:bound_hase},  and Corollary~\ref{thm:hase_cor}, 
we make the following assumptions.

\begin{assumption}\label{assume:HASE} \ 

\begin{enumerate}[(i)]

\item The sampling distribution $\zeta$ dominates the uniform distribution on $S$, that is, 
$$ P\left( \bar{Y}_0^{HASE} \leq y \right) \geq  P\left( \bar{Y}_0^{HAS1} \leq y \right)
{\rm \ for\ } y_* < y \leq y^*.$$
\label{assume:HASEinitial}
\item The bettering probability in HAS-E is bounded below by a positive constant, that is, for some positive $\gamma$,  $0<\gamma \leq1$, 
$$b(y) \geq \gamma {\rm \ for\ } y_* < y \leq y^*.$$
\label{assume:HASEbetteringprob}
\end{enumerate}
\end{assumption}


\vspace{0.5mm}
\noindent Theorem \ref{thm:hase_stoch_dom} proves  stochastic dominance of HAS-E over HAS1.

%
%
%


\begin{theorem}
\label{thm:hase_stoch_dom}
	Given the conditions in Assumption~\ref{assume:HASE} and setting $R_k=R$ for all $k$ as, 
	\begin{equation} \label{eq:second_rep_bound_hase}
	R = \left( \frac{ \sqrt[n]{q} \cdot 2 \cdot  \sigma \cdot z_{\alpha/2}  }{(1 - \sqrt[n]{q}) \cdot r_{y_* + \epsilon} \cdot \mathcal{K}_{q} } \right)^2  
	\end{equation}
	then
	 $\bar{Y}_k^{HASE}$ stochastically dominates $\bar{Y}_k^{HAS1}$, that is,
	$$P(\bar{Y}_k^{HASE} \leq y ) \geq P(\bar{Y}_k^{HAS1} \leq y ) \text{  for  } k = 0, 1, \ldots $$
	\noindent
	where $y_* < y \leq y^* $.
\end{theorem}

\noindent 
\textbf{Proof:} The proof is provided in Appendix~\ref{sec:proofs_hase}.

\vspace{2mm}
Since $\bar{Y}_k^{HASE}$ stochastically dominates $\bar{Y}_k^{HAS1}$, and the finite-time performance of HAS1 is captured in  \cite{Bulger1998,Wood2001}, we can bound the finite-time behavior of HAS-E.  

We are particularly interested in the expected behavior.
We  next derive an upper bound on the expected number of HAS-E iterations to achieve a sample point within a target level set $S_{y_*+\epsilon}$ for $\epsilon >0$, denoted $E[N^{HASE}_I(y_* +\epsilon)]$. The proof relies on the stochastic dominance of HAS-E over HAS1 in Theorem \ref{thm:hase_stoch_dom}, and uses an upper bound on HAS1 iterations,  as in \eqref{eq:has_bound_origin}. Theorem~\ref{thm:bound_hase} also expresses the expected number of function evaluations, including replications, needed to achieve a sample point within the target level set $S_{y_*+\epsilon}$, denoted $E[N^{HASE}_R(y_* +\epsilon)]$. 

%

\begin{theorem} \label{thm:bound_hase}
	An upper bound on the expected number of HAS-E iterations until reaching a value of $y_* + \epsilon$ or better, for $\epsilon >0$, is given by,
	\begin{equation} \label{eq:bound_has_update}
	\begin{split}
	& E[N^{HASE}_I(y_* +\epsilon)] \leq 1 + \left( \frac{1}{\gamma \cdot (1-\alpha) \cdot   q  }\right)  ln \left( \frac{\nu(S)}{\nu(S_{y_* + \epsilon})} \right)
	\end{split}
	\end{equation}
	
%
%

\noindent and an upper bound on the expected number of HAS-E function evaluations including replications is
\begin{equation} \label{eq:expectedN_R}
E[N^{HASE}_R(y_* +\epsilon)] \leq \left( \left(  \frac{q}{1-q}  + \frac{-ln(q) }{\left(1-q\right)^2} \cdot n \right) \left( \frac{2 \cdot  \sigma \cdot z_{\alpha/2} }{r_{y_* + \epsilon} \cdot \mathcal{K}_{q} }\right)  \right)^2 \cdot E[N^{HASE}_I(y_* +\epsilon)].
\end{equation}

\end{theorem}

\noindent
\textbf{Proof:} See the proof in Appendix~\ref{sec:proofs_hase}.

 \vspace{2mm}
 The expressions in \eqref{eq:bound_has_update} and \eqref{eq:expectedN_R} provide insight into the value of replications. Focusing on the impact of $q$ in  \eqref{eq:bound_has_update}, we see that the expected number of iterations decreases as $q$ increases, indicating some benefit to a large value of $q$ with a relatively tight upper confidence bound.   However, the expected number of replications in \eqref{eq:expectedN_R} indicates that the large number of replications needed for a large value of $q$ overshadows the benefit.  This can be interpreted as a tradeoff between sampling from a larger than needed level set (with loose upper confidence bound and fewer replications) and sampling from a more accurate estimate of the current level set (with tight upper confidence bound and more replications). This leads us to consider algorithms that use few replications as long as the estimation approaches the true function value as the algorithm approaches the global minimum. 
 
We also see the impact of dimension $n$ relative to the number of replications embedded in the sampling strategy using the upper confidence bound estimate. The expected number of replications in \eqref{eq:expectedN_R}  indicates that sampling on the level set associated with  the upper confidence bound as opposed to the true level set  increases the number of function evaluations quadratically in dimension as opposed to the number needed  if the function were able to be evaluated exactly. Also, the dimension $n$ magnifies the difference comparing values of $q$, reinforcing the intuition that fewer replications are better. 

The following corollary couples this with a bound on $\nu(S)/\nu(S_{y_* + \epsilon})$ in terms of  dimension $n$ for a class of problems satisfying a Lipschitz constant.

\begin{corollary} \label{thm:hase_cor}
	
	When $S$  in (P) is a convex feasible region in $n$ dimensions with a diameter $d$ and $f(x)$ satisfies the Lipschitz condition with Lipschitz constant at most $\mathcal{L} $, then the expected number of iterations for HAS-E to reach a value $y_*+\epsilon$, $\epsilon > 0$ , is bounded by,
	
	\begin{equation}\label{eq:has_bound_lk}
	E[N^{HASE}_I(y_* +\epsilon)] \leq  1 + \left( \frac{n  }{\gamma \cdot (1-\alpha) \cdot   q  }  \right)    ln \left( \frac{\mathcal{L} \cdot  d }{ \epsilon} \right)		
	\end{equation}  

%
	
	\noindent
	and the expected number of function evaluations (including replications) to achieve a value of $y_* + \epsilon$ or better
	is upper-bounded by a cubic function of domain dimension,
	\begin{equation*}
	\begin{split}
	& E \left[ N^{HASE}_R(y_* +\epsilon) \right] \\ 
	& \leq \left( \left(  \frac{q}{1-q}  + \frac{-ln(q) }{\left(1-q\right)^2} \cdot n \right) \left( \frac{2 \cdot  \sigma \cdot z_{\alpha/2} }{r_{y_* + \epsilon} \cdot \mathcal{K}_{q} }\right)  \right)^2  
	\left( 1 +     \left( \frac{n  }{\gamma \cdot (1-\alpha) \cdot  q }  \right)      ln \left( \frac{\mathcal{L} \cdot  d }{ \epsilon} \right) \right)   \\ 
	  	& \sim O\left(n^3 \right) . \\ 
	\end{split}
	\end{equation*}


\end{corollary}
\noindent
\textbf{Proof:} The expression in (\ref{eq:bound_has_update}), combined with the bounds on $\left( \nu(S) / \nu(S_{y_* + \epsilon}) \right)$ in \cite{Wood2001,Zabinsky2003,Zabinsky1992}, produce the linear number of iterations.  Coupling this with the expression in \eqref{eq:expectedN_R} provides an upper bound that is cubic in dimension.
\hspace{3mm} \qedsymbol

\vspace{2mm}
This corollary provides a bound on the expected number of function evaluations including replications needed to obtain a value of $y_* +\epsilon$ or less  that is  a \textit{polynomial} function of the dimension, holding all other parameters constant.   This result generally extends the finite-time results of the HAS framework for problems with estimation. In the next section, we examine a framework based on an adaptive search framework that samples from a series of nested quantile level sets.

\section{Quantile Adaptive Search with Estimation (QAS-E)}
\label{sec:qas}

We now define a Quantile Adaptive Search with Estimation (QAS-E) which conceptualizes an optimization algorithm that samples according to a probability distribution parameterized by quantile. QAS-E utilizes a series of sampling distributions defined by density function $\zeta_k$ associated with quantile $\delta_k$ on iteration $k$. 



The motivation for incorporating a quantile as a parameter in the sampling distribution is to provide a finite-time analysis to aid in the development of algorithms that use quantiles in their adaptive mechanics \cite{HoOrdinal2000,Ho:OrdinalBook,JiangHuPeng2022quantilebased,CrossEntropy,zabinsky2019PBnB}. 
QAS-E differs from HAS-E in its sampling distribution; instead of sampling $X_{k+1}$ from the normalized restriction of $\zeta$ on $S_{\bar{y}_k^{high}}$ as in HAS-E, QAS-E always samples from $S$ however, the distribution $\zeta_k$ depends on a quantile parameter $\delta_k$.  The intuition is that it is relatively easy to sample a point in a level set associated with a high quantile, but it is challenging to sample a point from a level set associated with a low quantile.  Instead of attempting to hit a target level set associated with a low quantile on the first iteration, QAS-E allows the quantile to be reduced iteratively, thereby modifying the parameterized sampling distribution. The hope is that small changes in the reduction of the quantile parameter will aid in implementation. 

We draw an analogy to the use of a temperature parameter in the Boltzmann distribution.  When the temperature is high, it is relatively easy to sample from the Boltzmann distribution, but when the temperature is low, it is difficult to efficiently generate a sample point.  The idea is that it is computationally easier to approximate a Boltzmann distribution with a small change in temperature, gradually reducing the temperature. The analysis of Annealing Adaptive Search  provided insight that led to the development of an adaptive cooling schedule for simulated annealing \cite{Shen2005,Shen2007}. 

The following analysis of QAS-E provides insight into the computational potential for algorithms that focus on sampling from level sets with quantile estimators. A general challenge with implementation is selecting the $\delta_k$-quantile values and associated sampling distributions $\zeta_k$ for which an adaptive algorithm has desirable performance.
We parameterize the sampling distribution by a quantile value, denoted $\zeta_k(\delta_k)$.  This is analogous to the way the Boltzmann distribution is parameterized by temperature.  In the following analysis, the Boltzmann distribution is  a possible family for QAS-E.
\color{black}

 There is a relationship between a quantile value $\delta$ and the associated objective function value $y_\delta$.
 For a quantile value, $0 < \delta <  1$, let the associated level set be denoted
\noindent
\begin{equation} S_{\delta}= \{x \in S: f(x) < y_\delta \} 
\label{eq:Sdelta}
\end{equation}
\noindent
where $y_{\delta}$ is the $\delta$-quantile of the domain $S$, or explicitly,
\begin{equation} y_{\delta} = \argmin_{y_* < y \leq y^* } P(f(X) \leq y) \geq \delta 
\label{eq:ydelta}
\end{equation}
when $X$ is uniformly sampled on $S$.  
\noindent 

When sampling according to the probability distribution $\zeta_k(\delta_k)$ on $S$, we relate the probability of landing inside of a level set $S_{\delta}$ to  the probability of achieving an objective function value of $y_{\delta}$ or better through an integral, as
$$ P \left( Y_k \leq y_{\delta} \right)=P(X_{k} \in S_{\delta})  = \int_{S_{\delta}} \zeta_k(\delta_k)(x) \cdot dx  $$
\noindent
where $X_k$ is drawn from $\zeta_k(\delta_k)$ and $Y_k=f(X_k)$.

\begin{figure}[!ht]
	\centering
	\includegraphics[height=3.1in]{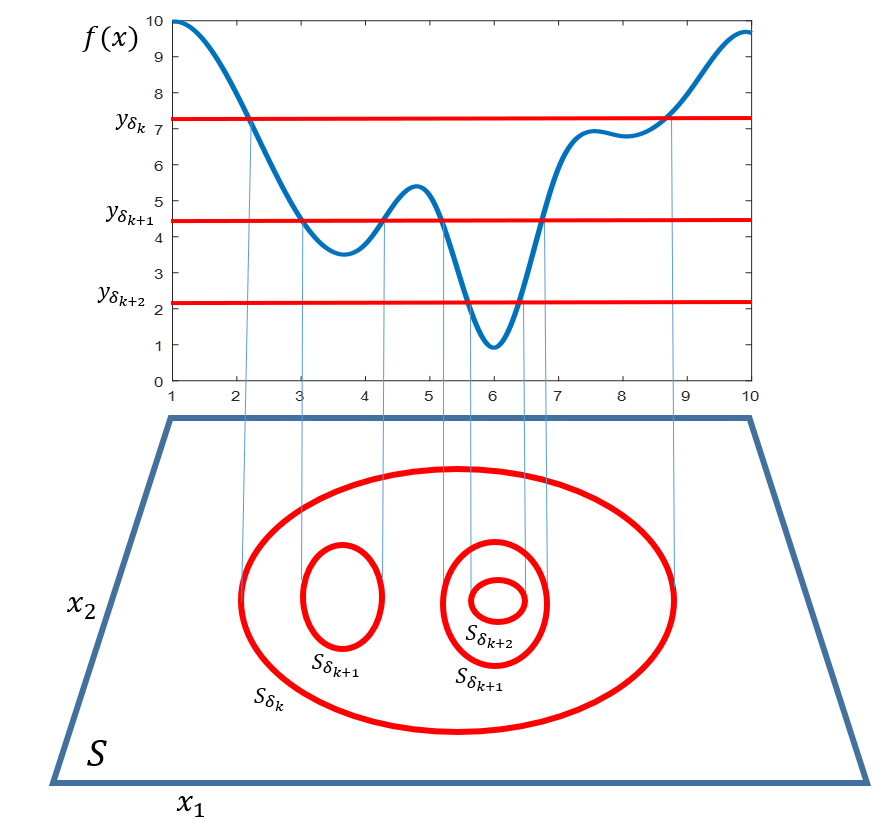}
	\caption{An illustration of three nested quantile level sets with $\delta_{k} > \delta_{k+1} > \delta_{k+2}$. Quantile Adaptive Search seeks to sample from the level set $S_{\delta_{k}}$ on the $k$th  iteration. } 
	\label{fig:Concept1st2}
\end{figure}

To illustrate the general form of QAS-E, see Figure \ref{fig:Concept1st2} with  three level sets associated with decreasing quantile levels $\delta_{k+2} < \delta_{k+1} < \delta_{k}$ so that $ S_{\delta_{k+2}} \subset S_{\delta_{k+1}} \subset S_{\delta_{k}} $. The sampling distribution $\zeta_k(\delta_k)$ is chosen to maintain some minimum probability of sampling within the associated level set $S_{\delta_{k}}$ at each iteration $k$. Therefore, the iterative selection of a quantile $\delta_{k}$ can be seen as a mechanism for focusing the sampling distribution on nested quantile level sets. 

QAS-E requires input parameters $\alpha$ and  $\sigma$, as in HAS-E, and a sequence of parameterized sampling distributions $\{\zeta_k(\delta_k), k=0,1,\ldots \}$.
It also requires a sequence of the number of replications on iteration $k$, i.e., $\{R_k,  k=0,1,\ldots \}$. The analysis of QAS-E requires several conditions on the sampling distributions stated in Assumption~\ref{assume:qas}, and discussed later.
We  formally write the algorithm based on the selection of parameterized sampling distributions  $\zeta_k(\delta_k)$.

\vspace{4mm}
\begin{samepage}    
\noindent
\textbf{Quantile Adaptive Search with Estimation (QAS-E)}
\begin{itemize} 
	\item \textbf{Step 0:} 
	Sample $X_0$ in $S$ according to the probability distribution 
	$\zeta_0(\delta_0)$ 
	on $S$. Conduct $R_0$ independent replications of the function at the initial selected point, i.e., $g(X_0, \chi_r) $ for $r = 1, \ldots, R_0$. Estimate the value $\hat{y}_0^{high} $ as in (\ref{eq:basic_estimate}) and set $\bar{y}_0^{high} = \hat{y}_0^{high}$. Set $\bar{Y}_0=f(X_0)$.	Set $k = 0$. 
	
	\item \textbf{Step 1:} Update the parameter quantile $\delta_{k+1}$ and its sampling probability distribution $\zeta_{k+1}(\delta_{k+1})$. 
	Generate $X_{k+1}$ from the probability distribution $\zeta_{k+1}(\delta_{k+1})$ on $S$.
Perform $R_k$ independent replications of $g(X_{k+1},\chi_r)$ for $r=1, \ldots, R_k$ (if  $X_{k+1} \neq  X_k$) and estimate $\hat{y}_{k+1}^{high}$ as in (\ref{eq:basic_estimate}). Then update  
$$	\bar{Y}_{k+1} =
\begin{cases}
f(X_{k+1}) & \text{if}\ f(X_{k+1}) < \bar{Y}_{k} \\
\bar{Y}_{k} & \text{otherwise}
\end{cases} $$
	and its associated upper confidence bound estimate,
\begin{equation}
	\bar{y}_{k+1}^{high}= \left\{ \begin{array}{ll}
	     \hat{y}_{k+1}^{high} & {\rm if}\  f(X_{k+1}) < \bar{Y}_{k} \\
	     \bar{y}_{k}^{high} & {\rm otherwise.}
	     \end{array} \right.
	     \nonumber
\end{equation}
	\item \textbf{Step 2:} If a stopping criterion is met, stop. Otherwise, increment $k$ and return to Step $1$.
\end{itemize}
\end{samepage}

\vspace{2mm}
Note that when there is no noise in the objective function, then no replications are needed and $\hat{y}_k^{high}$ can be replaced with the true function value $y_k $  in the algorithm.

The QAS-E algorithm iteratively samples from a sequence of distributions parameterized by a quantile value.   The intent is for the distributions to increase the chances of generating improving sets, much in the same way that the Boltzmann distribution with a temperature parameter increases its focus on improving level sets.  At each iteration, the upper confidence bound estimate $\bar{y}_k^{high}$ has an associated quantile value (denoted $\bar{\delta}_k^{high}$) through \eqref{eq:Sdelta} and \eqref{eq:ydelta}, that may help inform the choice of quantile parameter.  
The relationship between quantile as a parameter of the sampling distribution in QAS-E may guide the implementation of adaptive random search methods in an analogous way that AAS aided in developing a cooling schedule for simulated annealing. 
 

For the performance analysis of QAS-E, we  make the following assumptions regarding the sampling distribution with quantile parameter, $\zeta_k(\delta_k)$  on iteration $k$. The conditions in Assumption~\ref{assume:qas} ensure that each sampling distribution does no worse than the previous one at generating improving points.
 Assumptions~\ref{assume:qas}({\it i})-({\it ii}) are 
similar to Assumptions~1({\it i})-({\it ii})  for HAS-E.

\begin{assumption} \label{assume:qas} \ 

\begin{enumerate}[(i)]

\item \textit{The sampling distribution $\zeta_{k}(\delta_{k})$ dominates the uniform distribution, that is,} 
	\begin{align}
	&P \left( \bar{Y}_{k}^{QASE}  \leq y  \right) \geq \left( \bar{Y}_{0}^{HAS2}  \leq y  \right)
\end{align} 
\textit{for any iteration $k$ and $y_* < y   \leq y^*$. This requirement forces each sampling distribution to be more focused on improvement than the uniform distribution.  In effect this excludes distributions that are not able to sample from nested level sets better than uniform, perhaps due to local behavior.}	
\label{assume:qas-uniform}

	\item \textit{The probability of improving on the current upper bound estimate $\bar{y}_k^{high}$ when sampling from the probability distribution $\zeta_{k+1}(\delta_{k+1})$ is bounded below by some minimum probability $\gamma$}, 
	\begin{equation}
	P\left( \bar{Y}_{k+1} \leq \bar{y}_k^{high} |   \bar{Y}_{k} = \bar{y}_k \right) \geq \gamma  
	\end{equation}
where  $0<\gamma\leq 1$.	
	\textit{We require that the sampling distribution has a minimum probability of improvement. This requirement forces each updated sampling distribution to maintain some probability of sampling within the improving quantile level set associated with the upper bound estimate.}
		\label{assume:qas-better}

	\item \textit{The conditional probability that the distribution $\zeta_{k+1}(\delta_{k+1})$ samples within a lower level set given that the previous sampled value was $y_k$, is non-increasing in $\bar{y}_k$ for all $k$.  This condition can be written as}
	\begin{equation} 
	P \left( \bar{Y}_{k+1} \leq y |  \bar{Y}_{k}  = \bar{y}_k  \right)  
	\geq P \left( \bar{Y}_{k+1} \leq y |  \bar{Y}_{k}  = \bar{y}^\prime_k  \right)
	\end{equation}	
\textit{where $\bar{y}_k <  \bar{y}^\prime_k$. The sampling distributions cannot perform worse having observed a better point, e.g., $\bar{y}_k <  \bar{y}^\prime_k$. This prevents the sampling distribution from getting \enquote{stuck} at local minima by arriving at some small value that makes the sampling of further improvement almost impossible. }
	\label{assume:qas-nonincreasing}
\end{enumerate}
\end{assumption}

We now present an analysis of the performance of QAS-E that parallels that of HAS-E.  First, in Theorem~\ref{thm:qas_dom}, we show that the iterates of QAS-E stochastically dominate those of a special case of the standard HAS algorithm (called HAS2). Then, we use the special case HAS2 in Theorem~\ref{thm:qas_finitetime} to provide an upper bound on the expected number of QAS-E iterations and expected number of function evaluations including replications to achieve an optimal point with a function value within $\epsilon$ of the optimal value $y_*$.
%



%
The special case HAS2 uses uniform sampling, i.e., $\zeta^{HAS2} \sim {\rm Uniform}$, and the bettering probability is chosen to be constant for all $y$,
\begin{equation}
b^{HAS2}(y)= \gamma \cdot \left( 1 - \alpha \right) \cdot  q  
\label{eq:HAS2betteringprob}
\end{equation} 
where $0 < \gamma \leq 1$,
 $0<\alpha<1$,  and $0<q<1$.


Let $\bar{Y}_k^{QASE}$ be the best sampled value by QAS-E on the $k$th  iteration, 
and let $\bar{Y}_k^{HAS2}$ be the best sampled value on the $k$th iteration of HAS2. 
We show that QAS-E stochastically dominates HAS2 in Theorem~\ref{thm:qas_dom}.

%

\begin{theorem}\label{thm:qas_dom}
	Given the three conditions in Assumption~\ref{assume:qas} and setting $R_k=R$ for all $k$ as in (\ref{eq:second_rep_bound_hase}),  
	then $ \bar{Y}_k^{QASE}$ stochastically dominates  $ \bar{Y}_k^{HAS2}$, that is:
	$$P( \bar{Y}_k^{QASE} \leq y ) \geq P( \bar{Y}_k^{HAS2} \leq y ) \text{  for } k = 0, 1,  \ldots, $$
	\noindent 
	where $ y_* < y \leq y^* $.
\end{theorem}

\noindent
\textbf{Proof:} The proof is similar to the proof of Theorem~\ref{thm:hase_stoch_dom}, but is provided in Appendix \ref{sec:proofs_qas} for completeness. 
\vspace{2mm}

Theorem~\ref{thm:qas_finitetime} provides upper bounds on the expected number of QAS-E iterations and expected number of function evaluations including replications to achieve a point within a target level set $S_{y_* +\epsilon}$.  
Notice the bounds in Theorem~\ref{thm:qas_finitetime} are the same as in Theorem~\ref{thm:bound_hase}, suggesting the importance of an effective sampling distribution. 
	
\begin{theorem} \label{thm:qas_finitetime}
	An upper bound on the expected number of QAS-E iterations until the value of $y_* +\epsilon$ or better is sampled, for $\epsilon > 0$, is given by,

	\begin{equation}\label{eq:upper_bound_qas}
	\begin{split}
	& E[N_I^{QASE}(y_* + \epsilon)]
	\leq 1 + \left( \frac{1}{\gamma \cdot (1-\alpha) \cdot   q } \right)  ln \left( \frac{\nu(S)}{\nu(S_{y_* + \epsilon})} \right).\\
	\end{split}
	\end{equation}
	\noindent
	
\noindent and an upper bound on the expected number of QAS-E function evaluations including replications is

\begin{equation} \label{eq:QASEexpectedN_R}
E[N^{QASE}_R(y_* +\epsilon)] \leq \left( \left(  \frac{q}{1-q}  + \frac{-ln(q) }{\left(1-q\right)^2} \cdot n \right) \left( \frac{2 \cdot  \sigma \cdot z_{\alpha/2} }{r_{y_* + \epsilon} \cdot \mathcal{K}_{q} }\right)  \right)^2 \cdot E[N^{QASE}_I(y_* +\epsilon)].
\end{equation}

\end{theorem}

\noindent 
\textbf{Proof:} The proof is similar to that of Theorem~\ref{thm:bound_hase}.
%
 
 The final corollary is analogous to Corollary~\ref{thm:hase_cor}, and states that, when the problem (\textit{P}) satisfies certain conditions, the expected number of QAS-E iterations is bounded by a linear function in dimension, and the expected number of QAS-E function evaluations including replications is cubic in dimension.

\begin{corollary} \label{thm:qas_cor}
	
	When $S$  in ({\textit P}) is a convex feasible region in $n$ dimensions with a diameter $d$ and $f(x)$ satisfies the Lipschitz condition with Lipschitz constant at most $\mathcal{L} $, then the expected number of iterations for QAS-E to reach a value $y_*+\epsilon$, $\epsilon > 0$ , is bounded by,
	
	\begin{equation}\label{eq:qas_bound_lk}
	E[N_I^{QASE}(y_* + \epsilon)] 
	\leq  
	 1 + \left( \frac{n}{\gamma \cdot (1-\alpha) \cdot  q } \right)  ln \left( \frac{\mathcal{L} \cdot  d }{ \epsilon} \right)	
	\end{equation}  

\noindent and the expected number of function evaluations (including replications) to achieve a value of 	$y_*+\epsilon$ or better is bounded by a cubic function of domain dimension, 

	\begin{equation*}
	\begin{split}
	& E \left[ N^{QASE}_R(y_* +\epsilon) \right] \\ 
& \leq \left( \left(  \frac{q}{1-q}  + \frac{-ln(q) }{\left(1-q\right)^2} \cdot n \right) \left( \frac{2 \cdot  \sigma \cdot z_{\alpha/2} }{r_{y_* + \epsilon} \cdot \mathcal{K}_{q} }\right)  \right)^2  
	\left( 1 +    \left( \frac{n  }{\gamma \cdot (1-\alpha) \cdot   q }  \right)      ln \left( \frac{\mathcal{L} \cdot  d }{ \epsilon} \right) \right)   \\ 
	  	& \sim O\left(n^3 \right) . \\ 
	\end{split}
	\end{equation*}

\end{corollary}

\noindent
\textbf{Proof:}  The proof is similar to that of Corollary~\ref{thm:hase_cor}.
\hspace{3mm} \qedsymbol

\vspace{2mm}

\vspace{2mm}

The analysis of QAS-E parallels that for HAS-E, and highlights the result relating the estimation with a confidence bound and the performance related to the sampling distribution.  By making assumptions on the consistency of a sequence of sampling distributions, it is clear that there is flexibility in choosing a parameterized sampling distribution, however, the assumptions must be satisfied.   In this paper we emphasize using quantiles as parameters, however,  the Boltzmann distribution parameterized by temperature satisfies the assumptions too.  Thus, a version of AAS with estimation is captured in the analytical results.

%
%

\section{Discussion and Conclusion}
\label{sec:discussion}

%

We provide a framework for modeling adaptive random search when problems require estimation of the objective function. Hesitant Adaptive Search with Estimation  has a provable finite-time bound on the expected number of function evaluations until a specified $\epsilon$ above the minimum value is reached. Under certain conditions, the expected number of function evaluations including replications is bounded by a cubic function of dimension. 

Furthermore, we introduce a Quantile Adaptive Search with Estimation that extends HAS-E to an adaptive random search that always samples on $S$, but adapts the sampling distribution.  A difference between HAS-E and QAS-E is that HAS-E samples according to a normalized distribution restricted to nested level sets defined by the estimated function values, whereas QAS-E samples on the entire feasible region but parameterizes the sampling distribution to focus on nested level sets defined by quantiles.  QAS-E has similar finite-time results as HAS-E controlling the number of replications and iterations. 
The parameterized sampling distribution of QAS-E is analogous to use of the Boltzmann distribution with a temperature parameter in Annealing Adaptive Search. The analysis of QAS-E can be used to add estimation to AAS.

Future research will use the analysis of QAS-E to develop a means for adaptively setting quantile parameters in a sampling distribution, similar to how the analysis of AAS was used to derive an analytical cooling schedule for simulated annealing. The choice of sampling distributions must satisfy the assumptions put forth in Assumption~\ref{assume:qas}, however, the analysis may inspire algorithms like Probabilistic Branch and Bound \cite{zabinsky2019PBnB}, or other methods  
which either explicitly or implicitly attempt to sample from within quantiles of an objective function. It may be possible to use the analysis to inform reinforcement learning methods \cite{JiangHuPeng2022quantilebased}.

An insight that the analysis of HAS-E and QAS-E provides is that the value of consistent improvement in the sampling distribution is more important than the number of replications needed to achieve a close estimate of the objective function at points evaluated during the process. The bounds on expected function evaluations lead us to consider algorithms that use a few replications at the expense of sampling from a larger than needed level set. However, it is still important that the algorithm converges to the true global minimum. Future research will consider the Single Observation Stochastic Algorithm (SOSA) \cite{SOSA2018OR,SOSA2020ORletters} to combine estimation using a shrinking ball with the sampling distribution.

\clearpage
\begin{appendices}	\label{app:proofs}

%
%
%
%

\section{Proofs of Theorems for HAS-E Analysis}
\label{sec:proofs_hase}
\subsection{Proof of Theorem \ref{thm:lower_bound_convex}}

%

\noindent
\textbf{Proof of Theorem \ref{thm:lower_bound_convex}}:

\noindent 
For any value $y_k$ such that $y_* + \epsilon < y_k \leq y^*$, we start by defining an $n$-ball $\mathcal{B}_{y_k} $ as the largest $n$-ball centered at $x_*$ such that $\mathcal{B}_{y_k}  \subseteq S_{y_k}$ and let $r_{y_k}$ be its radius. We note that $0 < \nu(\mathcal{B}_{y_k}) \leq \nu(S_{y_k}) $. For any value $\hat{y}_{k}^{high}$, we define $\mathcal{B}_{\hat{y}_{k}^{high}}$ as the smallest $n$-ball centered at $x_*$ such that $S_{\hat{y}_{k}^{high} } \subseteq \mathcal{B}_{\hat{y}_{k}^{high} }$ and let $r_{\hat{y}_{k}^{high}}$ be the radius of  $\mathcal{B}_{\hat{y}_{k}^{high} }$. 

We  examine two cases. First, if  $\hat{y}_k^{high} - y_{k} \leq \kappa_{q}$ then $\frac{\nu(S_{y_{k}})}{\nu( S_{\hat{y}_k^{high}})} >  q $ by definition in \eqref{eq:K_q}, and the theorem is proved. 

\begin{figure}[!ht]
	\centering
		\includegraphics[height=3in]{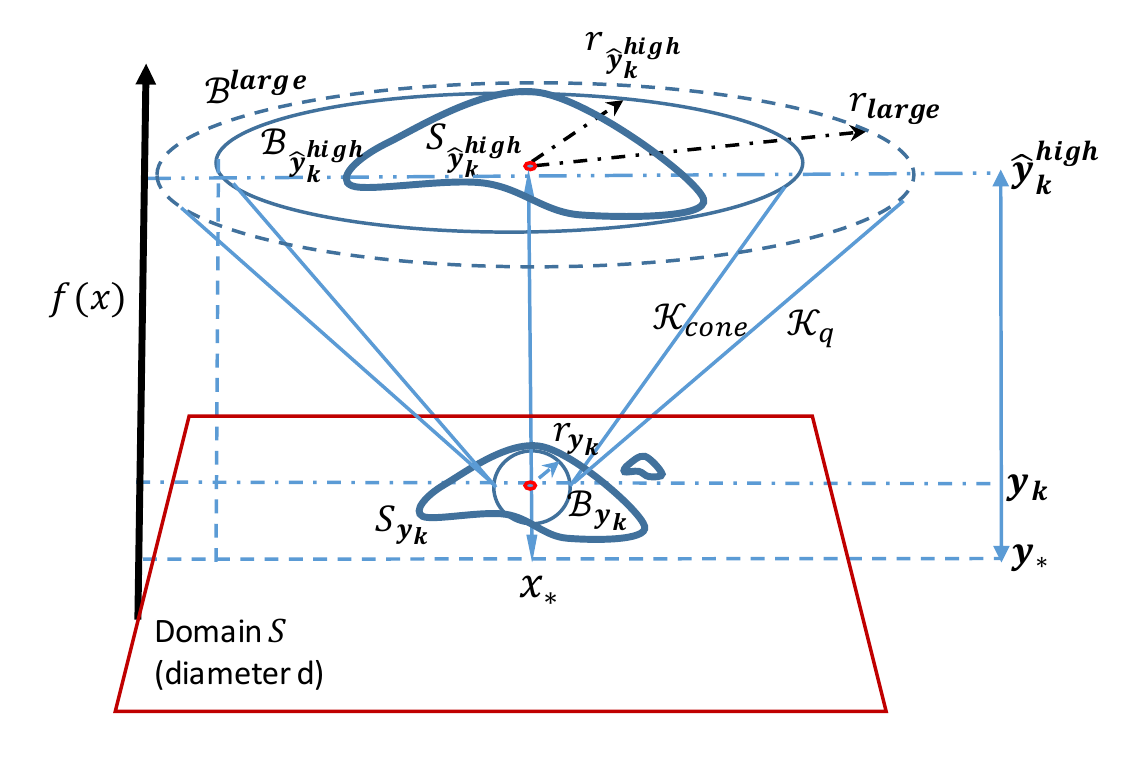}
	\caption{An illustration of the largest $n$-ball inscribed in $S_{y_k}$ $(\mathcal{B}_{y_k})$, the smallest $n$-ball inscribing $S_{\hat{y}_{k}^{high} }$ $(\mathcal{B}_{\hat{y}_{k}^{high}})$,  and a larger ball defined by the slope $\mathcal{K}_{q}$ ($\mathcal{B}^{large}$). } 	\label{fig:ConeProof}
\end{figure}

Second, consider  $\hat{y}_k^{high} - y_{k} >  \kappa_{q}$. 
We define $\mathcal{K}_{cone} = \frac{\hat{y}_{k}^{high}  - y_k}{  r_{\hat{y}_{k}^{high}} - r_{y_k} }$, which can be interpreted as the slope of the cone connecting the two balls, see Figure~\ref{fig:ConeProof}. We also  
write 
$$r_{\hat{y}_{k}^{high}} = r_{y_k} + (\hat{y}_{k}^{high}  - y_k) / \mathcal{K}_{cone}.  $$ 

Since $\hat{y}_k^{high} - y_{k} > \kappa_{q}$, the numerator of $\mathcal{K}_{cone}  $ is greater than the numerator of  $\mathcal{K}_{q}$ as in (\ref{eq:K_q}), and, since $d >r_{\hat{y}_{k}^{high}} - r_{y_k} $ by definition of the diameter, we have $\mathcal{K}_{cone}  > \mathcal{K}_{q}$. Note that $\mathcal{K}_{q}$ is independent of the value $y_k$.

We define $\mathcal{B}^{large} $ as an $n$-ball centered at $x_*$ with radius $r_{large} $, where $ r_{large} = r_{y_k} + (\hat{y}_{k}^{high} - y_k) / \mathcal{K}_{q}  $. Here we see that $r_{\hat{y}_{k}^{high}} \leq r_{large}  $ since $\mathcal{K}_{q} \leq \mathcal{K}_{cone}$. Therefore $S_{\hat{y}_{k}^{high} } \subset \mathcal{B}_{\hat{y}_{k}^{high} } \subset \mathcal{B}^{large} $, as illustrated in Figure \ref{fig:ConeProof}.

A lower bound on the ratios of volumes is constructed in terms of the dimension $n$, using  multi-dimensional geometry theorems \cite{kendall2004course},
$$ 	\frac{\nu(S_{y_{k}})}{\nu( S_{\hat{y}_k^{high}})} \geq    \frac{\nu(\mathcal{B}_{y_k})}{\nu(\mathcal{B}_{\hat{y}_{k}^{high} })} \geq \frac{\nu(\mathcal{B}_{y_k})}{\nu(\mathcal{B}^{large})} = \left( \frac{r_{y_k} }{r_{y_k} + \frac{\hat{y}_k^{high} - y_k}{\mathcal{K}_{q}}  }  \right)^n . $$
Since $ \hat{y}_k^{high} - y_k \leq \frac{ 2 \cdot \sigma \cdot z_{\alpha/2}}{\sqrt{R}}$, as given in the theorem statement, we have the following lower bound,
$$ \frac{\nu(S_{y_{k}})}{\nu( S_{\hat{y}_k^{high}})} \geq \left( \frac{r_{y_k}}{r_{y_k} + \frac{ 2 \cdot \sigma \cdot z_{\alpha/2}}{\mathcal{K}_{q} \sqrt{R}}  }  \right)^n . $$

We want to determine $R$, such that
$$ \left( \frac{r_{y_k}}{r_{y_k} + \frac{ 2 \cdot \sigma \cdot z_{\alpha/2}}{\mathcal{K}_{q} \sqrt{R}}  }  \right)^n \geq q.$$
\noindent
Taking a positive exponent of the two positive expressions, yields,
$$ \left( \frac{r_{y_k}}{r_{y_k} + \frac{ 2 \cdot \sigma \cdot z_{\alpha/2}}{\mathcal{K}_{q} \sqrt{R}}  }  \right) \geq \sqrt[n]{q}$$
\noindent 
and multiplying by a positive term,
$$ r_{y_k} \geq \sqrt[n]{q} \cdot \left( r_{y_k} + \frac{ 2 \cdot \sigma \cdot z_{\alpha/2}}{\mathcal{K}_{q} \sqrt{R}} \right)$$
\noindent 
and subtracting,
$$ r_{y_k} - \sqrt[n]{q} \cdot r_{y_k} \geq   \sqrt[n]{q} \cdot  \frac{ 2 \cdot \sigma \cdot z_{\alpha/2}}{\mathcal{K}_{q} \sqrt{R}} $$
\noindent 
and manipulating positive values,
$$ \sqrt{R} \geq  \sqrt[n]{q} \cdot  
\frac{ 2 \cdot \sigma \cdot z_{\alpha/2}}{ \left(  ( 1- \sqrt[n]{q})  \cdot r_{y_k} \cdot \mathcal{K}_{q}   \right) } $$
\noindent
and squaring both sides, yields
%
\begin{equation} 
R \geq \left( \frac{\sqrt[n]{q} \cdot  2  \cdot  \sigma  \cdot z_{\alpha/2} }{(1 - \sqrt[n]{q}) \cdot r_{y_k} \cdot \mathcal{K}_{q} } \right)^2. \nonumber
\end{equation}
Therefore, (\ref{eq:q_bound}) holds if $R \geq \left( \frac{ \sqrt[n]{q} \cdot  2 \cdot  \sigma \cdot z_{\alpha/2}  }{(1 - \sqrt[n]{q}) \cdot r_{y_k} \cdot \mathcal{K}_{q} } \right)^2$. Finally, since $y_k \geq y_*+\epsilon$, then $r_{y_k} \geq r_{y_*+\epsilon}$, and  (\ref{eq:q_bound}) holds if 
$$ R  \geq \left( \frac{2 \cdot \sqrt[n]{q} \cdot  2 \cdot  \sigma \cdot  z_{\alpha/2}  }{(1 - \sqrt[n]{q}) \cdot r_{y_*+\epsilon} \cdot \mathcal{K}_{q} } \right)^2$$
\noindent
which proves  Theorem \ref{thm:lower_bound_convex}. 
\hspace{3mm} \qedsymbol

\subsection{Additional Lemma}
Theorems \ref{thm:hase_stoch_dom} and \ref{thm:qas_dom} make use of Lemma $30$ from \cite{Shen2005}, which is repeated here for convenience.

\begin{lemma}\label{thm:base_lemma}
	\textbf{(cf. \cite{Shen2005})} \hspace{1mm} Let $ {\bar{Y}_k^A, k = 0, 1, 2, \ldots }$ and $ {\bar{Y}_k^B, k = 0, 1, 2,\ldots } $ be two sequences of objective function values generated by algorithms A and B respectively for solving a minimization problem, such that $\bar{Y}_{k+1}^A \leq \bar{Y}_k^A$ and $\bar{Y}_{k+1}^B \leq \bar{Y}_k^B$ for $k = 0, 1, \ldots $. For $y_* < y,z \leq y^*$ and $k = 0,1,  \ldots$, if 
	\begin{enumerate}
		\item $P(\bar{Y}_{k+1}^A \leq y | \bar{Y}_k^A =z  ) \geq P(\bar{Y}_{k+1}^B\leq y  | \bar{Y}_k^B =z  ) $
		\item $P(\bar{Y}_{k+1}^A \leq y | \bar{Y}_k^A =z  )$  is non-increasing in $z$, and 
		\item $ P(\bar{Y}_{0}^A \leq y ) \geq  P(\bar{Y}_{0}^B \leq y )$
	\end{enumerate}
	then $P(Y_{k}^A \leq y ) \geq P(Y_{k}^B \leq y)$ for $k=0,1,  \ldots$ and $y_* < y \leq y^*$. 
\end{lemma}

\noindent 
Proof of Lemma \ref{thm:base_lemma} can be found in \cite{Shen2005}.

\subsection{Proof of Theorem~\ref{thm:hase_stoch_dom}}


\noindent
\textbf{Proof of Theorem \ref{thm:hase_stoch_dom}: } Using the notation in HAS-E on the $k$th iteration, and based on Lemma \ref{thm:base_lemma}, as in \cite{Shen2005}, if the following conditions hold for $y_* < y, \bar{y}_k \leq y^*$ and $k = 0,1, \ldots $,
\begin{enumerate}[(I)]
	\item $P(\bar{Y}_{k+1}^{HASE}  \leq y | \bar{Y}_k^{HASE} =\bar{y}_k   ) \geq P(\bar{Y}_{k+1}^{HAS1}  \leq y | \bar{Y}_k^{HAS1} =\bar{y}_k   ) $
	\item $P(\bar{Y}_{k+1}^{HASE}  \leq y | \bar{Y}_k^{HASE} =\bar{y}_k   ) $  is non-increasing in $\bar{y}_k $, and 
	\item $ P(\bar{Y}_{0}^{HASE} \leq y) \geq  P(\bar{Y}_{0}^{HAS1}  \leq y )$
\end{enumerate}

\noindent 
then $P(\bar{Y}_k^{HASE} \leq y ) \geq P(\bar{Y}_k^{HAS1} \leq y ) \text{  for  } k = 0,1,  \ldots $  and for $y_*< y \leq y^*$. 
\vspace{2mm}

The first step is to prove (I). 
When $y \geq \bar{y}_k$, (I) is true trivially (since the conditional probability equals one on both sides). Now, when $y < \bar{y}_k$, we bound the left-hand side of the expression in (I), as,
\begin{equation}
\begin{split}
& P(\bar{Y}_{k+1}^{HASE} \leq y | \bar{Y}_{k}^{HASE}  = \bar{y}_k  ) \geq \gamma \cdot  P(\bar{Y}_{k+1}^{HASE} \leq y |  
 \bar{Y}_{k+1}^{HASE} \leq \bar{y}_k^{high},
\bar{Y}_{k}^{HASE} = \bar{y}_k  )
\end{split}
\label{eq:HASE_betters}
\end{equation}
%
%
\noindent
where we condition on the event that HASE \enquote{betters}, that is, that HASE samples from the normalized restriction of $\zeta$ on $S_{\bar{y}_k^{high}}$, and consequently $\bar{Y}_{k+1}^{HASE} \leq \bar{y}_k^{high}$,
which occurs with probability at least $\gamma$ by the bound on the bettering probability  in Assumption~\ref{assume:HASE}({\it \ref{assume:HASEbetteringprob}}).
 
%

We next consider the event $\{ \bar{y}_k \leq \bar{y}_k^{high} \}$, which occurs with probability at least $1-\alpha$ by  \eqref{eq:conf}. We rewrite \eqref{eq:HASE_betters} as,
\begin{equation} \label{eq:hase_breakdown}
\begin{split}
&P(\bar{Y}_{k+1}^{HASE} \leq y |  \bar{Y}_{k}^{HASE} = \bar{y}_k  )\\
& \hspace{4mm} \geq \gamma \cdot (1-\alpha) \cdot P(\bar{Y}_{k+1}^{HASE} \leq y |  \{ \bar{y}_k \leq \bar{y}_k^{high} \}, \bar{Y}_{k+1}^{HASE} \leq \bar{y}_k^{high})
\end{split}
\end{equation}
dropping the condition $\bar{Y}_{k}^{HASE} = \bar{y}_k$ because it is captured in the other conditions.
 
 From Assumption~\ref{assume:HASE}({\it \ref{assume:HASEinitial}}), we have 
\begin{equation}
\begin{split}
&P(\bar{Y}_{k+1}^{HASE} \leq y |  \bar{Y}_{k}^{HASE} = \bar{y}_k  )\\
& \hspace{4mm} \geq \gamma \cdot (1-\alpha) \cdot P(\bar{Y}_{0}^{HAS1} \leq y |  \{ \bar{y}_k \leq \bar{y}_k^{high} \}, \bar{Y}_{0}^{HAS1} \leq \bar{y}_k^{high})
\end{split}
\end{equation}
because $P(\bar{Y}_{0}^{HASE} \leq y) \geq P(\bar{Y}_{0}^{HAS1} \leq y)$ implies $P(\bar{Y}_{k+1}^{HASE} \leq y)\cdot  P(\bar{Y}_{k+1}^{HASE} \leq \bar{y}_k^{high})/P(\bar{Y}_{k+1}^{HASE} \leq \bar{y}_k^{high})\geq P(\bar{Y}_{0}^{HAS1} \leq y)\cdot P(\bar{Y}_{0}^{HAS1} \leq \bar{y}_k^{high})/P(\bar{Y}_{0}^{HAS1} \leq \bar{y}_k^{high})$.  Therefore, we can write

\begin{align} 
P(\bar{Y}_{k+1}^{HASE} \leq y |  \bar{Y}_{k}^{HASE} = \bar{y}_k  )
&\geq \gamma \cdot (1-\alpha) \cdot P(\bar{Y}_{0}^{HAS1} \leq y |  \bar{Y}_{0}^{HAS1} \leq \bar{y}_k^{high},  \{ \bar{y}_k \leq \bar{y}_k^{high} \})  \nonumber \\
&= \gamma \cdot (1-\alpha) \cdot \frac{P(\bar{Y}_{0}^{HAS1} \leq y)}{P(\bar{Y}_{0}^{HAS1} \leq \bar{y}_k^{high})}  \nonumber \\
&= \gamma \cdot (1-\alpha) \cdot \frac{P(\bar{Y}_{0}^{HAS1} \leq y)}{P(\bar{Y}_{0}^{HAS1} \leq \bar{y}_k)}  \cdot \frac{P(\bar{Y}_{0}^{HAS1} \leq \bar{y}_k)}{P(\bar{Y}_{0}^{HAS1} \leq \bar{y}_k^{high})} \nonumber \\
%
%
&\geq \gamma \cdot (1-\alpha)  \cdot \frac{\nu(S_{y})}{\nu(S_{\bar{y}_k}  )}  \cdot \frac{\nu(S_{\bar{y}_k})}{\nu(S_{\bar{y}_k^{high}}  )} \nonumber \\
&\geq \gamma \cdot (1-\alpha) \cdot  \frac{\nu(S_{y})}{\nu(S_{\bar{y}_k}  )} \cdot q .%
\label{eq:hase_breakdown_last}
\end{align}  
The last inequality makes use of  the lower bound developed in Theorem \ref{thm:lower_bound_convex}, $  \frac{\nu(S_{\bar{y}_k})}{ \nu(S_{\bar{y}_k^{high}}  ) }  \geq q $.\\

We similarly expand the expression for HAS1 in the right-hand side of (I), noting that HAS1 either improves or stays where it is, yielding,
 \begin{equation*}
P(\bar{Y}_{k+1}^{HAS1} \leq y | \bar{Y}_{k}^{HAS1} = \bar{y}_k  )
 = b^{HAS1}(\bar{y}_k)  P \left(  \bar{Y}_{k+1}^{HAS1} \leq y  |    \bar{Y}_{k}^{HAS1} = \bar{y}_k  \right) 
\end{equation*}
where ${X}_{k+1}^{HAS1}$ is sampled according to the normalized restriction of the uniform distribution on the improving level set.  Combining this with the
 bettering probability of HAS1,  $b(y)=\gamma \cdot (1-\alpha) \cdot q$, and when HAS1 \enquote{betters}, we have,
\begin{align}
&P(\bar{Y}_{k+1}^{HAS1} \leq y | \bar{Y}_{k}^{HAS1} = \bar{y}_k   )  
=   
\gamma \cdot (1-\alpha) \cdot  q \cdot \frac{\nu(S_{y})}{\nu(S_{\bar{y}_k})} . 
\label{eq:has1_breakdown_last}
 \end{align}

Combining \eqref{eq:hase_breakdown_last} and \eqref{eq:has1_breakdown_last}
proves condition (I).


\vspace{5mm}
We go on to prove (II), that $ P(\bar{Y}_{k+1}^{HASE} \leq y | \bar{Y}_k^{HASE} = \bar{y}_k  ) $ is non-increasing in $\bar{y}_k$. Suppose that $\bar{y}_k$ and $\bar{y}_k'$ are such that $\bar{y}_k < \bar{y}_k'$. To show (II) we want to show that: 
$$P(\bar{Y}_{k+1}^{HASE} \leq y | \bar{Y}_{k}^{HASE} = \bar{y}_k) \geq  P(\bar{Y}_{k+1}^{HASE} \leq y | \bar{Y}_{k}^{HASE} = \bar{y}_k'). $$
\noindent
The approach is to condition on the value of  $\bar{y}_k^{high}$, and since HAS-E samples on $S_{\bar{y}_k^{high}}$ in Step 2 of the algorithm, we know that $P(\bar{Y}_{k+1}^{HASE} \leq y  | \bar{Y}_k^{HASE} = \bar{y}_k, \bar{y}_k^{high} = u   )$  is non-increasing, therefore, we have,
\begin{align}
&P(\bar{Y}_{k+1}^{HASE} \leq y | \bar{Y}_{k}^{HASE} = \bar{y}_k) \hspace{80mm} \nonumber \\
&\hspace{5mm} = \int_{ - \infty}^{\infty} P(\bar{Y}_{k+1}^{HASE} \leq y | \bar{Y}_{k}^{HASE} = \bar{y}_k, \bar{y}_k^{high} = z ) \cdot dP(\bar{y}_k^{high} \leq z  | \bar{Y}_{k}^{HASE} = \bar{y}_k )  \nonumber
\end{align}
and because $\int_{-\infty}^{z}  dP(\bar{Y}_{k+1}^{HASE} \leq y | \bar{Y}_{k}^{HASE} = \bar{y}_k, \bar{y}_k^{high} = u ) = P(\bar{Y}_{k+1}^{HASE} \leq y | \bar{Y}_{k}^{HASE} = \bar{y}_k, \bar{y}_k^{high} = z ) - P(\bar{Y}_{k+1}^{HASE} \leq y | \bar{Y}_{k}^{HASE} = \bar{y}_k, \bar{y}_k^{high} = - \infty ) $, and since $P(\bar{Y}_{k+1}^{HASE} \leq y | \bar{Y}_{k}^{HASE} = \bar{y}_k, \bar{y}_k^{high} = -\infty ) = 1$ (trivially), we substitute $ P \left( \bar{Y}_{k+1}^{HASE} \leq y | \bar{Y}_{k}^{HASE} = \bar{y}_k, \bar{y}_k^{high} = z  \right) $ as follows, 
\begin{align}
&\hspace{7mm} = \int_{-\infty}^{\infty} \left( 1 + \int_{-\infty}^{z}  dP(\bar{Y}_{k+1}^{HASE} \leq y | \bar{Y}_{k}^{HASE} = \bar{y}_k, \bar{y}_k^{high} = u )\right) \cdot dP(\bar{y}_k^{high} \leq z  | \bar{Y}_{k}^{HASE} = \bar{y}_k ) \nonumber 
\end{align}
and reversing the order of integration, we get
\begin{align}
&\hspace{5mm} = 1 +  \int_{-\infty}^{\infty} \int_{u}^{\infty} dP(\bar{y}_k^{high} \leq z  | \bar{Y}_{k}^{HASE} = \bar{y}_k ) \cdot dP(\bar{Y}_{k+1}^{HASE} \leq y | \bar{Y}_{k}^{HASE} = \bar{y}_k, \bar{y}_k^{high} = u ) \nonumber
\end{align}
$$ \hspace{8mm}  = 1 +  \int_{-\infty}^{\infty} (1 - P(\bar{y}_k^{high} \leq u  | \bar{Y}_{k}^{HASE} = \bar{y}_k ) ) \cdot dP(\bar{Y}_{k+1}^{HASE} \leq y | \bar{Y}_{k}^{HASE} = \bar{y}_k, \bar{y}_k^{high} = u )  \hspace{20mm}   $$ 
however, since $dP(\bar{Y}_{k+1}^{HASE} \leq y | \bar{Y}_{k}^{HASE} = \bar{y}_k, \bar{y}_k^{high} = u ) \leq  0$, and since $P(\bar{Y}_{k+1}^{HASE} \leq y | \bar{Y}_{k}^{HASE} = \bar{y}_k, \bar{y}_k^{high} = u )$ is non-increasing in $\bar{y}_k^{high}$, and since, $P(\bar{y}_k^{high} \leq u  | \bar{Y}_{k}^{HASE} = \bar{y}_k ) \geq P(\bar{y}_k^{high} \leq u  | \bar{Y}_{k}^{HASE} = \bar{y}_k' )$, 
\noindent 
the probability that $\bar{y}_k^{high} $ is lower than $u$ is always greater for $\bar{y}_k < \bar{y}_k'$, then 
$$ \hspace{10mm} \geq  1 +  \int_{-\infty}^{\infty} (1 - P(\bar{y}_k^{high} \leq u  | \bar{Y}_{k}^{HASE} = \bar{y}_k' ) ) \cdot dP(\bar{Y}_{k+1}^{HASE} \leq y | \bar{Y}_{k}^{HASE} = \bar{y}_k, \bar{y}_k^{high} = u ) \hspace{20mm} $$
which is equivalent to
$$ \hspace{10mm} = 1 +  \int_{-\infty}^{\infty} \int_{u}^{\infty} (dP(\bar{y}_k^{high} \leq z  | \bar{Y}_{k}^{HASE} = \bar{y}_k' ) ) \cdot dP(\bar{Y}_{k+1}^{HASE} \leq y | \bar{Y}_{k}^{HASE} = \bar{y}_k, \bar{y}_k^{high} = u )  \hspace{95mm} $$
\noindent
and reversing the order of integration:
$$ \hspace{9mm} = 1 +  \int_{-\infty}^{\infty} \int_{-\infty}^{z}  dP(\bar{Y}_{k+1}^{HASE} \leq y | \bar{Y}_{k}^{HASE} = \bar{y}_k, \bar{y}_k^{high} = u )  \cdot (dP(\bar{y}_k^{high} \leq z  | \bar{Y}_{k}^{HASE} = \bar{y}_k' ) ) \hspace{65mm}  $$
$$ \hspace{8mm} =  \int_{-\infty}^{\infty} P(\bar{Y}_{k+1}^{HASE} \leq y | \bar{Y}_{k}^{HASE} = \bar{y}_k, \bar{y}_k^{high} = z )  \cdot dP(\bar{y}_k^{high} \leq z  | \bar{Y}_{k}^{HASE} = \bar{y}_k' )  \hspace{40mm}   $$
therefore, since $  P(\bar{Y}_{k+1}^{HASE} \leq y | \bar{Y}_{k}^{HASE} = \bar{y}_k, \bar{y}_k^{high} = z ) =  P(\bar{Y}_{k+1}^{HASE} \leq y | \bar{Y}_{k}^{HASE} = \bar{y}_k', \bar{y}_k^{high} = z ) $ , we write:
$$ \hspace{8mm} \geq  \int_{-\infty}^{\infty} P(\bar{Y}_{k+1}^{HASE} \leq y | \bar{Y}_{k}^{HASE} = \bar{y}_k', \bar{y}_k^{high} = z )  \cdot dP(\bar{y}_k^{high} \leq z | \bar{Y}_{k}^{HASE} = \bar{y}_k' ) \hspace{31mm} $$
$$  \hspace{8mm} =  P(\bar{Y}_{k+1}^{HASE} \leq y | \bar{Y}_{k}^{HASE} = \bar{y}_k') \hspace{96mm} $$
\noindent
which proves (II).

\vspace{2mm}

Lastly, condition (III) from Lemma \ref{thm:base_lemma} is true by Assumption~\ref{assume:HASE}({\it \ref{assume:HASEinitial}}) that  $P(\bar{Y}_0^{HASE}\leq y) \geq P(\bar{Y}_0^{HAS1} \leq y)$.
This proves the theorem through reference to Lemma \ref{thm:base_lemma}. \qquad
\qedsymbol

%

\subsection{Proof of a Second Lemma}

\noindent The following lemma is used in the proof of 
Theorem~\ref{thm:bound_hase}.

\begin{lemma}\label{thm:ratio_inequality}
	For a given constant $a$ such that $0 < a < 1$, and a variable $n \geq 1$, then the function $f(n)= \frac{a^{1/n}}{1-a^{1/n}}$ is bounded by a linear function of $n$, that is,
	\begin{equation}\label{eq:main_ineq} f(n) =  \frac{a^{1/n}}{1-a^{1/n}} \leq \frac{a}{1-a}  + \frac{-ln(a)  }{\left(1-a\right)^2} \cdot n.   \end{equation}
\end{lemma}

\noindent
\textbf{Proof of Lemma \ref{thm:ratio_inequality}:}
%
%
This bound is developed by proving that the derivative of $ \frac{a^{1/n}}{1-a^{1/n}}$ w.r.t $n$ obtains a finite maximum value over the range of $n \in [1, \infty]$. Using this upper bound on the derivative, a linear function is determined to bound the expression. First, note that when $n \geq 1$, the function is continuous and the derivative is defined. We take the first derivative of $f(n)$, yielding 
\begin{equation}\label{eq:first_deriv} \frac{d f(n)}{d n} = \frac{ -\frac{a^{1/n} ln(a)}{n^2} \cdot (1 -a^{1/n}) - \frac{a^{1/n} ln(a)}{n^2} \cdot (a^{1/n}) }{ (1 - a^{1/n})^2} = - \frac{a^{1/n} \cdot ln(a)}{n^2 \cdot( 1-a^{1/n} )^2} \end{equation} 
\noindent 
which is positive when $n > 0$. We go on to find a maximum value. 

Examine part of the denominator in (\ref{eq:first_deriv}), let 
$ \mathbf{d}(n) = n \cdot( 1-a^{1/n} )$
and we see that from Halley's Theorem:
$$\lim_{n \rightarrow \infty}  n \cdot( 1-a^{1/n} ) = - ln(a) $$
and at $n=1$, then 
$$\mathbf{d}(1)  = (1-a). $$
We first examine the first derivative of $\mathbf{d}(n)$,
$$ \mathbf{d}'(n) =  \frac{d \mathbf{d}(n)  }{ d n} = -a \cdot ^{1/n} + \frac{a^{1/n}\cdot ln(a)}{n} + 1 $$ 
with $ \lim_{n \rightarrow \infty} \mathbf{d}'(n)  =  0  $ and at $n=1$, then 
\begin{equation}\label{eq:a_expresion} \mathbf{d}' (1) = -a + a \cdot ln(a) + 1 = a\cdot(ln(a) - 1 )  + 1  . \end{equation} 
\noindent
Note that $\mathbf{d}'(1) > 0$ for $\forall a \in (-\infty, \infty)$ since (\ref{eq:a_expresion}) reaches a minimum in $a$ of $0$ at $a =1$. 

\vspace{2mm}
\noindent
Next, we examine the second derivative of $\mathbf{d}(n)$,
$$ \mathbf{d}''(n) =  \frac{d^2 \mathbf{d}(n)  }{ dn^2} =  - \frac{a^{1/n}\cdot ln^2 (a)}{n^3} .$$
We note that $\mathbf{d}''(n) < 0$ for $0 < a < 1$ and $n > 1$. Since $0 > \mathbf{d}''(n)  $ then $\mathbf{d}'(n) > 0 $ is always positive since $\mathbf{d}' $ monotonically decreases from $a\cdot(ln(a) - 1 )  + 1 $ to $0$ as $n$ increases. Similarly, since $\mathbf{d}'(n) \geq 0 $ for $n > 1$, then $\mathbf{d}(n)$ monotonically increases for $n > 1$. Therefore $\mathbf{d}(n)$ obtains a \textit{minimum} at $n=1$ therefore
\begin{equation} \label{eq:first_ineq}
n\cdot( 1-a^{1/n} ) \geq (1-a) .
\end{equation}

\noindent
Returning to (\ref{eq:first_deriv}), we develop an upper bound since $-ln(a)$ is positive and $a^{1/n} < 1$ for $0 < a <1$,
$$  f'(n) =  - \frac{a^{1/n} \cdot ln(a)}{n^2 \cdot( 1-a^{1/n} )^2}  \leq  \frac{-ln(a)}{\left( n\cdot( 1-a^{1/n} ) \right)^2}  $$
\noindent
then using (\ref{eq:first_ineq}),
$$ f'(n)  \leq \frac{-ln(a) }{\left(1-a\right)^2} . $$

\noindent
Using an upper bound of the derivative, $f'(n)$, and the value of $f(1) = \frac{a}{1-a}$, an upper bound is determined for for $f(n)$ when $n \geq  1$ as  
$$ \frac{a^{1/n}}{1-a^{1/n}} \leq \frac{a}{1-a}  + \frac{-ln(a) }{\left(1-a\right)^2} \cdot n . $$ 
which completes the proof. \hspace{3mm}
\qedsymbol 

\subsection{Proof of Theorem~\ref{thm:bound_hase}}

\noindent
\textbf{Proof of Theorem \ref{thm:bound_hase}:}
By stochastic dominance in Theorem \ref{thm:hase_stoch_dom}, the expected number of iterations to achieve a value within $S_{y_* + \epsilon}$ for HAS-E is less than or equal to the number for HAS1. Since the bettering probability for HAS1 is $b(y) = \gamma \cdot ( 1- \alpha) \cdot  q$  for all $ y_* < y \leq y^* $, using (\ref{eq:has_bound_origin}), we have
\begin{align*}
&E[N^{HASE}_I(y_* +\epsilon)] \leq 1 + \int_{y_* + \epsilon}^{\infty} \frac{d\rho (t)}{\gamma \cdot ( 1 - \alpha) \cdot  q   \cdot p(t)}\end{align*}
and since HAS1 uses uniform sampling, i.e., $ p(y) = \frac{\nu(S_y)}{\nu(S)}$, we have 
\begin{align*}
& \hspace{28mm} = 1 + \frac{1}{\gamma \cdot (1-\alpha) \cdot   q  } \cdot ln \left( \frac{\nu(S)}{\nu(S_{y_* + \epsilon})} \right). 
\end{align*}
Using a constant number of replications $R$ for each iteration, yields
$$E[N^{HASE}_R(y_* +\epsilon)] = R \cdot E[N^{HASE}_I(y_* +\epsilon)].$$

\noindent 
To bound
$ E \left[ N^{HASE}_R(y_* +\epsilon) \right]$, we apply  
 the inequality in \eqref{eq:main_ineq} (in Lemma~\ref{thm:ratio_inequality} in Appendix~A)
 to get
 $$\frac{\sqrt[n]{q}}{(1-\sqrt[n]{q})} \leq 
  \frac{q}{1-q}  + \frac{-ln(q) }{\left(1-q\right)^2} \cdot n $$
 and combining this  \eqref{eq:second_rep_bound_hase}, 
we obtain
	\begin{equation*}
	\begin{split}
	& E \left[ N^{HASE}_R(y_* +\epsilon) \right] \\ 
	& \leq    \left( \left(  \frac{q}{1-q}  + \frac{-ln(q) }{\left(1-q\right)^2} \cdot n \right) \left( \frac{2 \cdot  \sigma \cdot  z_{\alpha/2}  }{r_{y_* + \epsilon} \cdot \mathcal{K}_{q} }\right)  \right)^2 \cdot
	E[N^{HASE}_I(y_* +\epsilon)].
	\end{split}
	\end{equation*}

This proves the theorem. \qquad \qedsymbol

\section{Proofs of Theorems for QAS-E Analysis}
\label{sec:proofs_qas}

%


The proofs of the theorems for QAS-E are similar to the proofs for HAS-E.  The QAS-E proofs are provided for completeness.

\vspace{2mm}

\noindent 
\textbf{Proof of Theorem \ref{thm:qas_dom}}

Similar to the proof of Theorem \ref{thm:hase_stoch_dom}, if the three conditions listed in Lemma \ref{thm:base_lemma} hold, 
%
%
then $P(\bar{Y}_k^{QASE} \leq y ) \geq P(\bar{Y}_k^{HAS2} \leq y ) \text{  for  } k = 0,1, \ldots $  and for $y_*< y \leq y^*$. 
\vspace{2mm}

We start by proving the first condition in Lemma \ref{thm:base_lemma}, that is, we show that 
$$P(\bar{Y}_{k+1}^{QASE}  \leq y | \bar{Y}_k^{QASE} =\bar{y}_k   ) \geq P(\bar{Y}_{k+1}^{HAS2}  \leq y | \bar{Y}_k^{HAS2} =\bar{y}_k   ) $$
for  $y_* < y, \bar{y}_k \leq y^*$ and $k = 0,1, \ldots $.
When $y \geq \bar{y}_k$,   $P(\bar{Y}_{k+1}^{QASE} \leq y | \bar{Y}_k^{QASE} =\bar{y}_k ) = P(\bar{Y}_{k+1}^{HAS2} \leq y | \bar{Y}_k^{QASE} =\bar{y}_k  ) =1 $, and the first condition holds. 

Now, when $y < \bar{y}_k$,
we bound the left-hand side of the expression in (I) by conditioning on the event 
that $X_{k+1}^{QASE} $ \enquote{betters}, that is, the event that $\bar{Y}_{k+1}^{QASE} \leq \bar{y}_k^{high}$, yielding
\begin{align}
P(\bar{Y}_{k+1}^{QASE} \leq y | \bar{Y}_k^{QASE} =\bar{y}_k ) 
     & \geq \gamma \cdot P(\bar{Y}_{k+1}^{QASE} \leq y | \bar{Y}_{k+1}^{QASE} \leq \bar{y}_k^{high}, \bar{Y}_k^{QASE} =\bar{y}_k ) 
     \label{eq:QASbound1}
 \end{align}
 by Assumption~\ref{assume:qas}({\it \ref{assume:qas-better}}).
 
 We next consider the event $\{\bar{y}_k \leq \bar{y}_k^{high}\}$, which occurs with probability at least $1-\alpha$, by
 \eqref{eq:conf}. We rewrite \eqref{eq:QASbound1} as,
 \begin{equation}
 \begin{split}
&P(\bar{Y}_{k+1}^{QASE} \leq y |  \bar{Y}_{k}^{QASE} = \bar{y}_k  )\\
& \hspace{4mm} \geq \gamma \cdot (1-\alpha) \cdot P(\bar{Y}_{k+1}^{QASE} \leq y | \bar{Y}_{k+1}^{QASE} \leq \bar{y}_k^{high}, \{\bar{y}_k \leq \bar{y}_k^{high}\}, \bar{Y}_k^{QASE} =\bar{y}_k ). 
\end{split}
\end{equation}
 From Assumption~\ref{assume:qas}({\it \ref{assume:qas-uniform}}), we have 
\begin{align} 
P(\bar{Y}_{k+1}^{QASE} \leq y |  \bar{Y}_{k}^{QASE} = \bar{y}_k  )
&\geq \gamma \cdot (1-\alpha) \cdot P(\bar{Y}_{0}^{HAS2} \leq y |  \bar{Y}_{0}^{HAS2} \leq \bar{y}_k^{high}).  \nonumber 
%
\end{align}

Making use of  the lower bound developed in Theorem \ref{thm:lower_bound_convex}, $  \frac{\nu(S_{\bar{y}_k})}{ \nu(S_{\bar{y}_k^{high}}  ) }  \geq q $, we have,

\begin{align} 
P(\bar{Y}_{k+1}^{QASE} \leq y |  \bar{Y}_{k}^{QASE} = \bar{y}_k  )
&\geq \gamma \cdot (1-\alpha) \cdot P(\bar{Y}_{0}^{HAS2} \leq y |  \bar{Y}_{0}^{HAS2} \leq \bar{y}_k^{high})  \nonumber \\
&= \gamma \cdot (1-\alpha) \cdot \frac{P(\bar{Y}_{0}^{HAS2} \leq y)}{P(\bar{Y}_{0}^{HAS2} \leq \bar{y}_k^{high})}  \nonumber \\
&= \gamma \cdot (1-\alpha) \cdot \frac{P(\bar{Y}_{0}^{HAS2} \leq y)}{P(\bar{Y}_{0}^{HAS2} \leq \bar{y}_k)}  \cdot \frac{P(\bar{Y}_{0}^{HAS2} \leq \bar{y}_k)}{P(\bar{Y}_{0}^{HAS2} \leq \bar{y}_k^{high})} \nonumber \\
&\geq \gamma \cdot (1-\alpha) \cdot 
\frac{\nu(S_{y})}{\nu(S_{\bar{y}_k}  )}  \cdot 
\frac{\nu(S_{\bar{y}_k})}{\nu(S_{\bar{y}_k^{high}}  )} \nonumber \\
&\geq \gamma \cdot (1-\alpha) \cdot 
\frac{\nu(S_{y})}{\nu(S_{\bar{y}_k}  )} \cdot q .%
\label{eq:qase_breakdown_last}
\end{align}

We similarly expand the expression for HAS2 in the right-hand side of (I), noting that HAS2 either improves or stays where it is, yielding,
 \begin{equation*}
P(\bar{Y}_{k+1}^{HAS2} \leq y | \bar{Y}_{k}^{HAS2} = \bar{y}_k  )
= b^{HAS2}(\bar{y}_k)  P \left(  \bar{Y}_{k+1}^{HAS2} \leq y  |   \bar{Y}_{k}^{HAS2} = \bar{y}_k   \right)
\end{equation*}
and since the bettering probability of HAS2 equals $\gamma \cdot (1-\alpha) \cdot 
q$, and when HAS2 \enquote{betters}, it samples uniformly on the improving level set, we have,
\begin{align}
&P(\bar{Y}_{k+1}^{HAS2} \leq y | \bar{Y}_{k}^{HAS2} = \bar{y}_k   )  
=   
\gamma \cdot (1-\alpha) \cdot 
q \cdot \frac{\nu(S_{y})}{\nu(S_{\bar{y}_k})} . 
\label{eq:has2_breakdown_last}
 \end{align}

Combining \eqref{eq:qase_breakdown_last} and \eqref{eq:has2_breakdown_last}
proves condition (I).

\noindent

The second condition in Lemma \ref{thm:base_lemma} is satisfied directly by Assumption~\ref{assume:qas}({\it \ref{assume:qas-nonincreasing}}). The third condition in Lemma \ref{thm:base_lemma} is satisfied  by Assumption~\ref{assume:qas}({\it \ref{assume:qas-uniform}}).  
This proves the theorem by Lemma \ref{thm:base_lemma}. \qedsymbol
\color{black}

\end{appendices}


\bibliographystyle{plain}
\bibliography{bibUpdateZZ}

\begin{thebibliography}{10}

\bibitem{BAS2007}
William~P. Baritompa, David~W. Bulger, and Graham~R. Wood.
\newblock Generating functions and the performance of backtracking adaptive
  search.
\newblock {\em Journal of Global Optimization}, 37(2):159--175, 2007.

\bibitem{boender1995stochastic}
C.~Guus~E. Boender and H.~Edwin Romeijn.
\newblock Stochastic methods.
\newblock In {\em Handbook of global optimization}, pages 829--869. Springer,
  1995.

\bibitem{Bulger1998}
David~W. Bulger and Graham~R. Wood.
\newblock {Hesitant adaptive search for global optimisation}.
\newblock {\em Mathematical Programming}, 81(1):89--102, 1998.

\bibitem{FuHandbook}
Michael~C. Fu.
\newblock {\em Handbook of Simulation Optimization}, volume 216.
\newblock Springer New York, New York, NY, 2015.

\bibitem{HoOrdinal2000}
Y.~C. Ho, C.~G. Cassandras, C.~H. Chen, and L.~Dai.
\newblock Ordinal optimisation and simulation.
\newblock {\em Journal of the Operational Research Society}, 51:490--500, 2000.

\bibitem{Ho:OrdinalBook}
Y.~C. Ho, Q.~C. Zhao, and Q.~S. Jia.
\newblock {\em Ordinal optimization: Soft optimization for hard problems}.
\newblock Springer, Berlin, Germany, 2007.

\bibitem{HuJiaqiao2012ASoS}
Jiaqiao Hu, Yongqiang Wang, Enlu Zhou, Michael~C. Fu, and Steven~I. Marcus.
\newblock A survey of some model-based methods for global optimization.
\newblock In {\em Optimization, Control, and Applications of Stochastic
  Systems}, pages 157--179. Birkh\"{a}user Boston, 2012.

\bibitem{JiangHuPeng2022quantilebased}
Jinyang Jiang, Jiaqiao Hu, and Yijie Peng.
\newblock Quantile-based policy optimization for reinforcement learning, 2022.
\newblock available on arXiv 2201.11463.

\bibitem{kendall2004course}
Maurice~G. Kendall.
\newblock {\em A Course in the Geometry of n Dimensions}.
\newblock Courier Corporation, 2004.

\bibitem{SOSA2018OR}
S.~Kiatsupaibul, R.~L. Smith, and Z.~B. Zabinsky.
\newblock Single observation adaptive search for continuous simulation.
\newblock {\em Operations Research}, 66:1713 -- 1727, 2018.

\bibitem{SOSA2020ORletters}
S.~Kiatsupaibul, R.~L. Smith, and Z.~B. Zabinsky.
\newblock Single observation adaptive search for discrete and continuous
  simulation.
\newblock {\em Operations Research Letters}, 48:666 -- 673, 2020.

\bibitem{locatelli2013global}
Marco Locatelli and Fabio Schoen.
\newblock {\em Global optimization: theory, algorithms, and applications},
  volume~15.
\newblock SIAM, 2013.

\bibitem{LocatelliSchoen2021}
Marco Locatelli and Fabio Schoen.
\newblock ({G}lobal) optimization: Historical notes and recent developments.
\newblock {\em EURO Journal on Computational Optimization}, 9:100012, 2021.

\bibitem{pardalos2000recent}
Panos~M. Pardalos, H.~Edwin Romeijn, and Hoang Tuy.
\newblock Recent developments and trends in global optimization.
\newblock {\em Journal of Computational and Applied Mathematics},
  124(1-2):209--228, 2000.

\bibitem{raphael2003direct}
Benny Raphael and Ian F.~C. Smith.
\newblock A direct stochastic algorithm for global search.
\newblock {\em Applied Mathematics and Computation}, 146(2-3):729--758, 2003.

\bibitem{RaphaelSmith}
Benny Raphael and Ian F.~C. Smith.
\newblock Global search through sampling using a {PDF}.
\newblock In {\em Stochastic Algorithms: Foundations And Applications}, volume
  2827, pages 71--82. Springer, 2003.

\bibitem{romeijn1994simulated}
H.~Edwin Romeijn and Robert~L. Smith.
\newblock Simulated annealing and adaptive search in global optimization.
\newblock {\em Probability in the Engineering and Informational Sciences},
  8(4):571--590, 1994.

\bibitem{CrossEntropy}
Reuven~Y. Rubinstein and Dirk~P. Kroese.
\newblock {\em The Cross-Entropy Method: A Unified Approach to Combinatorial
  Optimization, Monte-Carlo Simulation and Machine Learning}.
\newblock Springer, Cambridge, UK, 2004.

\bibitem{Shen2005}
Yanfang Shen.
\newblock {Annealing Adaptive Search With Hit-and-Run Sampling Methods for
  Global Optimization}.
\newblock {\em University of Washington Dissertation}, 2005.

\bibitem{Shen2007}
Yanfang Shen, Seksan Kiatsupaibul, Zelda~B. Zabinsky, and Robert~L. Smith.
\newblock {An analytically derived cooling schedule for simulated annealing}.
\newblock {\em Journal of Global Optimization}, 38(3):333--365, 2007.

\bibitem{BAS2006}
Graham~R. Wood, David~W Bulger, William~P. Baritompa, and D.~Alexander.
\newblock Backtracking adaptive search: Distribution of number of iterations to
  convergence.
\newblock {\em Journal of Optimization Theory and Applications},
  128(3):547--562, 2006.

\bibitem{Wood2001}
Graham~R. Wood, Zelda~B. Zabinsky, and Birna~P. Kristinsdottir.
\newblock Hesitant adaptive search: the distribution of the number of
  iterations to convergence.
\newblock {\em Mathematical Programming}, 89(3):479--486, 2001.

\bibitem{Zabinsky2003}
Zelda~B. Zabinsky.
\newblock {\em {Stochastic adaptive search for global optimization}}.
\newblock Kluwer Academic Publishers originally, Springer Science {\&} Business
  Media, 2003.

\bibitem{Zabinsky2010}
Zelda~B. Zabinsky, David Bulger, and Charoenchai Khompatraporn.
\newblock {Stopping and restarting strategy for stochastic sequential search in
  global optimization}.
\newblock {\em Journal of Global Optimization}, 46:273--286, 2010.

\bibitem{zabinsky2019PBnB}
Zelda~B. Zabinsky and Hao Huang.
\newblock A partition-based optimization approach for level set approximation:
  Probabilistic branch and bound.
\newblock In Alice Smith, editor, {\em Women in Industrial and Systems
  Engineering: Key Advances and Perspectives on Emerging Topics}. Springer
  Nature, 2020.

\bibitem{Zabinsky1992}
Zelda~B. Zabinsky and Robert~L. Smith.
\newblock {Pure adaptive search in global optimization}.
\newblock {\em Mathematical Programming}, 53(1-3):323--338, 1992.

\bibitem{Zabinsky1995}
Zelda~B. Zabinsky, Graham~R. Wood, Mike~A. Steel, and William~P. Baritompa.
\newblock Pure adaptive search for finite global optimization.
\newblock {\em Mathematical Programming}, 69(1-3):443--448, 1995.

\end{thebibliography}


\end{document}